\documentclass[11pt]{article}
\usepackage[utf8]{inputenc}
\usepackage[main=english]{babel}
\usepackage{mathtools,amssymb,amsthm,amsfonts,upgreek,dsfont,mathabx}
\usepackage{geometry}
\usepackage{enumerate}
\usepackage{amsmath}
\usepackage{url}
\usepackage{fancybox}
\usepackage{fullpage}
\usepackage[colorlinks, linktocpage, citecolor=blue, linkcolor=blue]{hyperref}
\usepackage{multicol} 
\usepackage{array}
\usepackage{xargs}
\usepackage{float}
\usepackage[prependcaption]{todonotes}
\usepackage[capitalise]{cleveref}
\usepackage{bbm}
\usepackage{subcaption}
\usepackage{comment}
\usepackage{todonotes}

\newcommandx{\Max}[2][1=]{\todo[inline, author={Maxime}, linecolor=green,backgroundcolor=green!25,bordercolor=green,#1]{#2}}

\newtheorem{theorem}{Theorem}[section]

\newtheorem{prop}{Proposition}[section]
\newtheorem{lem}{Lemma}[section]
\newtheorem{cor}{Corollary}[section]
\newtheorem{rem}{Remark}[section]
\newtheorem{assumption}{Assumption}
\newtheorem*{st_assumption}{Standing assumption}
\renewcommand{\P}{\mathbb{P}}
\newcommand{\R}{\mathbb{R}}
\newcommand{\N}{\mathbb{N}}
\newcommand{\E}{\mathbb{E}}

\newcommand{\Vnq}{V_n^{k,R,q}}
\newcommand{\Vnz}{V_n^{k,R,0}}

\DeclareMathOperator{\Var}{Var}
\DeclareMathOperator{\ESS}{ESS}
\DeclareMathOperator{\argmin}{argmin}

\definecolor{darkgreen}{rgb}{0.0, 0.5, 0.0}

\newcommand{\vertiii}[1]{{\left\vert\kern-0.25ex\left\vert\kern-0.25ex\left\vert #1 
    \right\vert\kern-0.25ex\right\vert\kern-0.25ex\right\vert}}

\newcommand{\logl}{\Phi}
\newcommand{\prior}{\pi_0}
\newcommand{\post}{\pi_n}
\newcommand{\logpostd}{I_n}

\newcommand{\lap}{\mathcal{L}_{\post}^K}

\newcommand{\impw}{w_n}
\newcommand{\unimpw}{w^*_n}

\title{Multi-modal Laplace based sampling}
 \author{
 {
 \large Maxime \textsc{Egéa}}\footnote{Email: \texttt{maxime.egea@sorbonne-universite.fr}, Sorbonne Universite} 
 \and
 {\large Hanyue \textsc{Gu}}\footnote{Email: \texttt{hanyue.gu@math.tu-freiberg}, TU Bergakademie Freiberg} 
 \and
 {\large Daniel \textsc{Rudolf}}\footnote{Email: \texttt{daniel.rudolf@uni-passau.de}, University of Passau} 
 \and
 {\large Bj\"orn \textsc{Sprungk}}\footnote{Email: \texttt{bjoern.sprungk@math.tu-freiberg.de}, TU Bergakademie Freiberg}
 }

\date{\today}
\begin{document}


\maketitle


\begin{abstract}
We study the multi-modal Laplace approximation, given by a weighted Gaussian mixture centred at the
local minima of the negative log-density, for sampling from probability measures that concentrate
around finitely many modes as an inverse-temperature parameter $n$ grows. 
Such distributions occur in  Bayesian inference in the small noise or large data regime.
We derive explicit, non-asymptotic bounds on the Hellinger distance between the multi-modal Laplace
approximation and the concentrating posterior. The bound decays like $n^{-1/2}$ with $n$ and the constants depend only polynomially on the dimension $d$.
Moreover, we analyse the Laplace mixture as a proposal distribution for self-normalised importance sampling and independent Metropolis--Hastings (IMH) targeting concentration multimodal distributions. 
We derive asymptotic results, as $n\to\infty$, on the effective sample size in case of importance sampling as well as on the average acceptance probability, lag-one autocorrelation, and spectral gap in case of IMH.
Numerical experiments illustrate and verify our theoretical findings.
\end{abstract}


\section{Introduction}


The extraction of information from complex, multi-modal distributions that are specified only implicitly through evaluations of their densities is a long-standing problem in  applied sciences. A promising approach is approximate sampling, that is, generating realizations from a probability measure that is `sufficiently close' to the target distribution. This methodology is frequently used in computational statistics, machine learning, and statistical physics, where the desired distribution often depends on an inverse-temperature parameter. As the inverse temperature increases, intuitively, the distribution becomes increasingly concentrated. While this concentration can make it possible to identify a particularly high probability area, it also poses significant challenges for sampling and can substantially reduce the efficiency of the sampling scheme. We propose a multi-modal Laplace mixture approximation that can serve as a proxy for the target distribution and under suitable assumptions remains capable to capture the asymptotic behavior when the `concentration' increases.

For formalizing the setting let $(\pi_n)_{n\in\N}$ be the sequence of
(increasingly) concentrating target posterior measures, each defined on $\R^d$ equipped with its Borel $\sigma$-algebra, of the form
\begin{equation}\label{eq:definitionofmu}
    \post(\mathrm{d}x) \coloneq \frac{1}{Z_n} \mathrm{e}^{-n\logl_n(x)} \pi_0(\mathrm{d}x), 
    \quad\text{with}\quad 
    Z_n \coloneq  \int_{\R^d} \mathrm{e}^{-n\logl_n(x)}\prior(\mathrm{d}x),
\end{equation}
where $\logl_n : \R^d\rightarrow [0,\infty)$ represents a (scaled) negative log-likelihood or potential 
which may depend on the `inverse temperature' $n\in \N$ and $\prior$ denotes a prior or reference probability measure on $\R^d$. 
As motivation of the setting we provide frameworks where such target distributions occur, for instance, as posterior measures in Bayesian inverse problems with real-valued observation vector $y$ of
\[
    Y = \mathcal{G}_n(X) + \eta_n \quad \text{with}\quad X \sim \pi_0.
\]
`Concentration behavior' for increasing $n$ of the corresponding posterior distribution is present in the following general scenarios:
\begin{enumerate}[a)]
\item \label{en: small_noise}
Assume that the observation operator $\mathcal{G}_n = G$ (with $o\in\N$, $y\in\R^o$ and $G\colon \R^d \to \R^o$ not depending on $n$) and vanishing noise $\eta_n \sim \mathcal{N}(0, n^{-1}\Gamma)$  (for positive definite symmetric matrix $\Gamma\in \R^{o\times o}$) determines
\[
    \logl_n(x) := \frac{1}{2}\|y - G(x)\|^2_{\Gamma^{-1}}, 
    \quad \text{with} \quad \|z\|_{\Gamma^{-1}} := \|\Gamma^{-1/2}z \|
\]
and Euclidean norm $\|z\|$ for $z\in\R^o$.
\item \label{en: large_data}
For $G$ given as in \ref{en: small_noise}), $y=(y_1,\dots,y_n)\in\R^{n\cdot o}$ for $y_1,\dots,y_n\in\R^o$, 
$\eta_n \sim \mathcal{N}(0, \sigma^2 \mathrm{Id}_{no})$ (with $\sigma>0$ and identity matrix $\mathrm{Id}_{no}\in \R^{n\cdot o \times n\cdot o}$) consider an observation operator of the form 
\[
\mathcal{G}_n = \underbrace{(G, \ldots, G)}_{n \text{-times}}\colon \R^d \to\R^{n\cdot o},\]
such that (for $i=1,\dots,n$ interpreting $y_i$ as realization of $G(X)+\xi_i$ with $\eta_n = (\xi_1,\dots,\xi_n)$) we have
\[
    \logl_n(x) := \frac{1}{2n \sigma^2} \sum_{k=1}^n \|y_k - G(x)\|^2.
\]
\end{enumerate}
 Scenario \ref{en: small_noise}) is referred to as small noise limit case and scenario \ref{en: large_data}) as large data limit case. Let us mention here that 
we will mostly focus on the small noise limit framework in the theoretical analysis but also consider the large data limit in our numerical experiments. 

For large $n\in\N$, a standard tool in Bayesian analysis, cf. \cite{gelman1995bayesian}, and Bayesian machine learning, cf. \cite{DaxbergerEtAl2021}, for dealing with a concentrated posterior measure $\post$ is to use their \emph{Laplace approximation (LA) $\mathcal L_{\post}$} as proxy or surrogate.  The LA is a Gaussian approximation to $\post$ centered around the (assumed) global maximum of its Lebesgue density. It has gained quite some attention in recent years. Particularly, there have been a number of results regarding explicit error bounds and convergence analysis for increasing $n\to\infty$, see for example \cite{ art:hasenpflug2023wasserstein, 
art:HelinKretschmannNonAs,  katsevich2024, katsevich2025,Spokoiny2023,art:SchillingBjornPhilipConvergLaplaceApprox, Spokoiny2025}. 

A particularly remarkable recent result obtained by Katsevich and Spokoiny in \cite{katsevich2025} provides an explicit upper bound of the
total variation distance between $\post$ and $\mathcal L_{\post}$ of the form $d \, n^{-1/2}$, which shows linear dependence on the dimension.
Notably, most of the previously cited papers consider the case  of uni-modal $\post$. In this work we analyze a suitable LA $\lap$ in the case of multiple but finitely many (here $K$) global minimizers of the potential $\Phi_n$. 
The, to our knowledge, first asymptotic result on this \emph{multi-modal case} for $\logl_n = \Phi$ being independent of $n$ can be found in \cite{art:Hwang_LaplaceMethod}. There, the author shows that when $\mathrm{argmin}_{x\in \R^d} \Phi(x) = \left\{x_1, \dots, x_K \right\}$ for $K \in \N$, the weak limit of $\post$ is a convex mixture of the Dirac measures $\delta_{x_1}, \dots, \delta_{x_K}$. Via the construction of explicit couplings, an extension of the former result is achieved in \cite{art:hasenpflug2023wasserstein}, where the $p$-Wasserstein distances (with $p\in[1,\infty)$) between $\post$ to the aforementioned Dirac mixture are estimated. In particular, the authors show a convergence of the form $n^{-1/2p}$ under suitable assumptions.\\[-0.75ex]

\textbf{Our contribution.} We study the behavior and benefits of the LA for the computation of expectations w.r.t.~$\post$ of quantities of interest by the use of Monte Carlo methods. 
In particular, we consider a Gaussian mixture LA for multimodal targets $\post$, cf. \cite{gelman1995bayesian}. In the first part of the paper, we generalize the convergence result of \cite[Theorem 2]{art:SchillingBjornPhilipConvergLaplaceApprox} to the case of multimodal LA and multimodal $\post$. 
Moreover, following ideas of \cite{art:HelinKretschmannNonAs} we also derive a non-asymptotic bound in Hellinger distance. 
A particular attention is given to explicitly stating all the constants involved in this bound, especially their dependence on the dimension $d$.
The second part of the paper is devoted to analyze LA-based methods for Monte Carlo integration w.r.t~$\post$. We study the use of the LA as proposal distribution for both self-normalized importance sampling and the independent Metropolis–Hastings algorithm and  analyze how their efficiency behaves as $n\to \infty$. For LA-based importance sampling, we establish convergence in probability of order $n^{-1+\delta}$, for arbitrary $\delta>0$, for the effective sample size tending to the total sample size as $n\to\infty$. For the independent Metropolis–Hastings algorithm we analyze the average acceptance probability, the lag-one autocorrelation, as well as the spectral gap. We provide, under suitable assumptions, a `convergence to 1' statement of the average acceptance rate and spectral gap, accompanied by a `convergence to 0' statement of the lag-one autocorrelation of order $n^{-1/2}$.
These theoretical findings are illustrated in various numerical examples.\\[-0.75ex]

\textbf{Outline.} In Section \ref{sec:setting} we introduce the general setting, assumptions and the multimodal LA. Section \ref{sec:conv} establishes the convergence of LA to $\pi_n$ in the Hellinger distance. In Section \ref{sec:LAIS} we analyze the LA-based importance sampling whereas in Section \ref{sec:LAIMH} we study the independent Metropolis--Hastings algorithm using the LA as proposal. Afterwards, both methods and statements are considered in several numerical experiments in Section \ref{sec:numerics}. Finally, Section \ref{sec:concl} concludes.


\section{Preliminaries and Laplace mixture approximation}\label{sec:setting}

\subsection{Notation and Preliminaries}
We denote the usual inner product on $\R^d$ and the induced Euclidean norm by $\langle \cdot, \cdot\rangle$ and $\|\cdot\|$, respectively. 
By $\vertiii{\cdot}$ we denote the matrix norm associated with the Euclidean norm and for a symmetric positive definite matrix $A \in \R^{d\times d}$ we define $\|\cdot\|_A= \|A^{1/2 } \cdot \|$.

For a given probability measure $\pi$ on $\R^d$ we denote by $L^p(\pi)$ the Lebesgue space of all $p$-integrable functions $f\colon \R^d \to \R$ and the corresponding $L^p$-norm by $\|\cdot\|_{L^p(\pi)}$.
For $f \in L^1(\pi)$ we set $\mathbb{E}_{\pi}[f] := \int_{\R^d} f \,\mathrm{d} \pi$ and for $f \in L^2(\pi)$ we denote $ \Var_{\pi}(f)
:= \mathbb{E}_{\pi}[(f - \mathbb{E}_{\pi}[f])^2]$.

In the following we work with three classical distances and divergences of probability measures $\pi$, $\widetilde \pi$ on $\mathbb R^d$, namely, the \emph{total variation distance $d_{\mathrm{TV}}$} and the  
\emph{Hellinger distance $d_{\mathrm{H}}$}, which are given by
\[
d_{\mathrm{TV}}(\pi,\widetilde \pi) 
\coloneqq 
\frac{1}{2} \int_{\mathbb R^d} 
\left|\frac{\mathrm d\pi}{\mathrm d\mu} - \frac{\mathrm d\widetilde \pi}{\mathrm d\mu}\right| \,\mathrm d\mu,
\qquad
d_{\mathrm H}(\pi, \widetilde \pi) 
\coloneqq 
\left(
\int_{\mathbb R^d} \left|\sqrt{\frac{\mathrm d\pi}{\mathrm d\mu}} - \sqrt{\frac{\mathrm d\widetilde \pi}{\mathrm d\mu}}\right|^2 \mathrm d\mu \right)^{1/2},
\]
where $\mu$ denotes a common dominating measure of $\pi$ and $\widetilde \pi$ on $\mathbb R^d$ 
and $\frac{\mathrm d\pi}{\mathrm d\mu}$, $\frac{\mathrm d\widetilde \pi}{\mathrm d\mu}$ the corresponding Radon--Nikodym derivatives, as well as the \emph{$\chi^2$-divergence} of $\widetilde \pi$ from $\pi$ defined as 
\[
    \chi^2(\widetilde \pi \parallel \pi)
    :=
    \int_{\R^d}
    \left(\frac{\mathrm d\widetilde \pi}{\mathrm d\pi} - 1\right)^2
    \ 
    \mathrm d\pi
    \in[0,\infty]
\]
for the case that $\widetilde \pi$ is absolutely continuous w.r.t.~$\pi$ and $\chi^2(\widetilde \pi \parallel \pi):=\infty$ otherwise.
Note that the $\chi^2$-divergence can also be infinite in the case when $\frac{\mathrm d\widetilde \pi}{\mathrm d\pi}$ exists and that it is not a metric in contrast to $d_{\mathrm{TV}}$ and $d_{\mathrm{H}}$. 
Regarding the latter it is well-known, cf. \cite[Equation (8)]{GibbSu2002}, that
\begin{equation}\label{eq:dTV_dH}
    \frac{d_{\mathrm H}^2(\pi, \widetilde \pi)}{2} 
    \leq 
    d_{\mathrm{TV}}(\pi,\widetilde{\pi}) 
    \leq 
    d_{\mathrm H}(\pi,\widetilde{\pi})
\end{equation}
and that the Hellinger distance allows to control differences of expectations of square integrable functions 
\begin{equation}\label{eq:pif_tildepif}
    \left|
    \mathbb{E}_{\pi}[f] - 
    \mathbb{E}_{\widetilde\pi}[f]
    \right| 
    \le \left( 2\|f\|_{L^2(\pi)}^2 + 2\|f\|_{L^2({\widetilde \pi})}^2\right)^{1/2} d_{\mathrm H}(\pi,\widetilde{\pi}),
    \qquad
    \forall 
    f \in L^2(\pi) \cap L^2(\widetilde{\pi}),
\end{equation}
see \cite[Lemma 21]{art:DashtiStuart_Bayesian}. 


\begin{rem}\label{rem:HellBound}
    Note, that the proof of the bound \eqref{eq:pif_tildepif} extends to general positive and finite measures $\pi,\widetilde \pi$ which might even have different total mass on $\R^d$. 
\end{rem}

Moreover, the Hellinger distance can be bounded by the $\chi^2$-divergence
\begin{equation}
    d_{\mathrm{H}}( \pi,\widetilde \pi) 
    \le 
    \chi^2 (\widetilde \pi \parallel \pi)
\end{equation}
see \cite[Section 2.4]{book:tsybakov_IntroToNonParamStat}. In general the converse inequality does not hold, but in case of $\left\|\frac{\mathrm d\widetilde \pi}{\mathrm d\pi}\right\|_\infty := \sup_{x\in\R^d} |\frac{\mathrm d\widetilde \pi}{\mathrm d\pi}(x)| < \infty$, we have
\begin{equation}\label{eq:chi2_hellinger}
    \chi^2 (\widetilde \pi \parallel \pi) 
    \le \left(\sqrt{\left\|\frac{\mathrm d\widetilde \pi}{\mathrm d\pi}\right\|_\infty}+1\right)^2 \ d_{\mathrm{H}}(\pi,\widetilde \pi)^2,
\end{equation}
which follows by using $(t-1)^2 = (\sqrt{t}-1)^2\ (\sqrt t+1 )^2$ for $t\geq0$ with $t = \frac{\mathrm d\widetilde \pi}{\mathrm d\pi}(x)$.

\subsection{Laplace mixture approximation}
We assume that the prior distribution $\prior$ in \eqref{eq:definitionofmu} is absolutely continuous with respect to the Lebesgue measure, with density $\prior : \R^d\rightarrow (0,\infty) $.
As a consequence, also $\post$ admits the representation in terms of a Lebesgue density
\[
    \post(\mathrm{d}x) \coloneq \frac{1}{Z_n} \mathrm{e}^{-n\logpostd(x)}\mathrm{d}x,
    \quad \text{with} \quad
    \logpostd \coloneqq \logl_n - \frac{1}{n} \log(\prior).
\]

As a standard assumption for the well-definedness of the Laplace approximation we impose the following differentiability requirements. For this let us denote with $\mathcal{C}^r\left(\R^d, [0,\infty) \right)$ the class of $r$-times continuously differentiable functions that map from $\R^d$ to $[0,\infty)$. The following condition is assumed to be satisfied throughout the whole paper.

\begin{st_assumption}
For each $n \in \N $, the function $\logpostd  \in \mathcal{C}^2\left(\R^d, [0,\infty) \right)$ has $K\in\N$ local minima $x_n^{(1)},\dots,x_n^{(K)} \in \R^d$. Moreover, for any $n\in\N$ and any $k=1,\ldots,K$ we set
\[
    \pi_n^{(k)} \coloneqq  \mathrm e^{-n I_n (x_n^{(k)} )} > 0,
    \qquad
    H_n^{(k)} 
    \coloneqq \nabla^2\logpostd (x_n^{(k)})
\]
and assume that the matrices $H_n^{(k)}$ are positive definite.
\end{st_assumption}

Under the \textbf{standing assumption}, we define the multi-modal version of the LA of $\pi_n$, with $n\in \N$, by a weighted Gaussian mixture (cf. \cite{gelman1995bayesian}), as
\begin{equation}\label{eq:Laplacemixture}
    \mathcal{L}_{\post}^K 
    \coloneqq 
    \sum_{k=1}^{K} 
    c_{n,k} \ 
    \mathcal{N}\left(x_n^{(k)}, \big(n H_n^{(k)}\big)^{-1} \right),
\end{equation}
with coefficients $c_{n,k} \in [0,1]$, satisfying $\sum_{k=1}^K c_{n,k}=1$, given by
\begin{equation*}
    c_{n,k} 
    \coloneqq \frac{\mathrm{det}\Big(\big(nH_n^{(k)}\big)^{-1}\Big)^{\frac{1}{2}} \pi_n^{(k)}}{\sum_{j=1}^K \mathrm{det}\left(\big(nH_n^{(j)}\big)^{-1}\right)^{\frac{1}{2}} \pi_n^{(j)}}.
\end{equation*}
This LA probability measure $\mathcal{L}_{\post}^K$ is absolutely continuous w.r.t. the Lebesgue measure and, for simplicity, we use the same notation for both the distribution and its probability density function $\lap$, that is, for any $x\in\R^d$ we have
\begin{align}
\label{al: Z_tilde}
    \lap(x)  & = \frac{1}{\widetilde Z_n^K}\sum_{k=1}^K \pi_{n}^{(k)}
    \exp\Big(-\frac{n}{2}\|x-x_n^{(k)}\|^2_{H_n^{(k)}} \Big), 
    \quad\text{with}\quad
      \widetilde Z_n^K  := \Big(\frac{2 \pi}{n} \Big)^{d/2}\sum_{j=1}^K  \det\big( \big(H_n^{(j)}\big)^{-1/2} \big) \pi_{n}^{(j)}.
\end{align}
\begin{rem}
We comment on the definition of the multi-modal LA. In principle, there are at least 
two ways to build it based on a Gaussian mixture since there are two types of local minima: those that vanish as $n$ goes to infinity, and those that become global in the limit. For instance, for the function $\logpostd\colon\R \to[0,\infty)$, defined by
\begin{equation*}
    \logpostd (x) = (x^2-1)^2 + \frac{1}{n}(x+1)^2,
\end{equation*}
we have $\arg \min_{x\in\R} I_n(x) =\Big \{ -1, \frac{1}{2}\Big(1+\sqrt{1-\frac{2}{n}}\Big)\Big\}$ and it can be shown that $\post$ converges in distribution toward the mixture of Diracs $\frac{\delta_{-1} + e^{-4}\delta_1}{1+e^{-4}} $. But the situation can also be different. Consider the only minorly changed
\begin{equation*}
    \logpostd (x) = (x^2-1)^2 + \frac{1}{\sqrt{n}}(x+1)^2, 
\end{equation*}
for which $\arg \min_{x\in\R} I_n(x) =\Big \{ -1, \frac{1}{2}\Big(1+\sqrt{1-\frac{2}{\sqrt{n}}}\Big) \Big\}$
and it can be shown that $\post$ converges in distribution toward $ \delta_{-1}$. This means that in the limit only a single minimum is present. As a consequence, to define a multi-modal LA in such a situation, we may either consider only the minima that become global in the limit or consider all local minima. Asymptotically, this does not change too much, since the coefficient corresponding to the Gaussian of the local minima that vanishes at infinity also tends toward $0$ as $n$ goes to infinity. However, taking all local minima into account seems more natural from a non-asymptotic point of view, since the approximation will be `closer' to the posterior.
\end{rem}

\section{Convergence in Hellinger distance}\label{sec:conv}

In this section, we investigate the convergence of the Hellinger distance between $\post$ and the LA $\lap$. For formulating the underlying assumptions we require for each local minimum $x_n^{(1)},\dots,x_n^{(K)}$ a neighborhood $\Vnq$, with $R,q\geq0$, defined as
\begin{equation*}
    \Vnq 
    \coloneqq
    \left\{ x\in \R^d : \|x-x_n^{(k)} \|^2_{H_n^{(k)}} \leq Rn^{-2q} \right\},
    \qquad
    n\in\N,\ k=1,\ldots,K.
\end{equation*}
Moreover, we introduce the `tail area' as
$
    \overline{V}_n^{R,q}
    \coloneqq 
    \R^d \setminus \left( \bigcup_{k=1}^K \Vnq \right).
$

\begin{rem}\label{rem:neighborhood}
From the definition it is visible that for any $k=1,\dots,K$ the sets $\Vnq$ `shrink' as $n^{-q}$ for increasing $n$. If the posterior $\post$ `concentrates rapidly around its local minima', then these neighborhoods `capture essentially all the mass of the posterior' as $n$ grows. 
\end{rem}
We are now in a position to impose differentiability conditions on $I_n$ that guarantee regular behavior in the neighborhoods.
\begin{assumption}\label{assum:Two}
    For any $n\in\N $ let $I_n \in \mathcal C^3(\R^d, [0,\infty))$ and for $q\geq0$ as well as for $R,\rho,M \in (0,\infty)$ assume that:
\begin{enumerate}[(i)]
    \item 
    \label{en: pairw_disj}
    For any $n\in\N$ the local minima neighborhoods $V_n^{1,R,q},\dots,V_n^{K,R,q}$ are pairwise disjoint. 
    \item 
    \label{en: unif_pos}
    The matrices $\big(H_n^{(k)}\big)_{k \in \{1, \dots , K \} }$ are $\rho$-uniformly positive definite, i.e., for all $k \in \{1, \dots , K \}$ and for each $n \in \N$ we have
    \begin{equation*}
        H_n^{(k)} \succeq \rho^2 \mathrm{Id}_d,
    \end{equation*}
    where the inequality is taken in the sense of positive definite matrices.
    \item \label{en: M_bounded}
    The third-order derivative tensors $D^3 I_n\left(x\right)$ are $M$-uniformly bounded on $\Vnq$, i.e., for all $k \in \{1, \dots , K \}$, each $n \in \N$, all $x\in \Vnq$, and all $h \in \R^d$ we have
    \begin{equation*}
        \left| D^3 \logpostd \left(x\right)(h,h,h)\right| \leq 6 M \|h\|^3.
    \end{equation*}
\end{enumerate}    
\end{assumption}

\begin{rem} 
\label{rem: norm_ineq}
Let us already here
note that \eqref{en: unif_pos} of Assumption~\ref{assum:Two}  implies that for each $n \in \N $, $k \in \{1, \dots, K \}$ and any $x\in\R^d$ holds 
\begin{equation*}
    \|x-x_n^{(k)} \| \le \vertiii{(H_n^{(k)})^{-1/2}}\|x-x_n^{(k)} \|_{H_n^{(k)}} \le \frac{1}{\rho} \|x-x_n^{(k)} \|_{H_n^{(k)}}. 
\end{equation*}
\end{rem}
\subsection{Asymptotic result for $q=0$}
Now, we establish a generalization of \cite[Theorem~2]{art:SchillingBjornPhilipConvergLaplaceApprox} to the multi-modal setting. We consider the case $q=0$, such that $\Vnq = \Vnz$ and note that the neighborhoods do not `shrink' w.r.t. $n$. Anyway, still we need to be able to control what happens in the `tail area', which motivates the following `domination condition'.


\begin{assumption} 
\label{assum:coerc}
For $R>0$ and $0<\varepsilon<1$, there is an $n^*\in\N$ such that for all $n\geq n^*$ we have
\begin{equation} \label{eq: add_req}
  \inf_{x\in \overline{V}_n^{R,0}} I_n(x)>\varepsilon,
  \qquad
     \sum_{j=1}^K \mathrm{det}\big((H_n^{(j)})^{-1}\big) \pi_n^{(j)}>\varepsilon,
\end{equation}
   and there is a function
   $g : \R^d \to [0,\infty)$, with $g^{1-\varepsilon}$ being Lebesgue integrable, such that for any $x \in \R^d$ and $n\geq n^*$ we have
    $
    \exp(-n\logpostd(x)) \leq g (x).
    $
\end{assumption}

\begin{rem}
If for a given $R>0$ we have 
\[
\lim_{n\to\infty} \inf_{x\in \overline{V}_n^{R,0}} I_n(x)>0
\quad \text{and} \quad
\lim_{n\to\infty} \sum_{j=1}^K \mathrm{det}\big((H_n^{(j)})^{-1}\big) \pi_n^{(j)} >0,
\]
then there is an $\varepsilon>0$ and $n^*\in\N$, such that \eqref{eq: add_req} holds for any $n\geq n^*$ as asked for in the previous assumption.
\end{rem}

 Note that the combination of Assumption~\ref{assum:Two} and Assumption~\ref{assum:coerc} is adapted from the requirements formulated in \cite{art:SchillingBjornPhilipConvergLaplaceApprox} to cover the multi-modal setting. Now we state one of our main theorems that is proven in Appendix~\ref{sec: proof_main_Helling}.
\begin{theorem} \label{theo:MainpropHelling2} 
    For $q=0$ and $\rho,M \in (0,\infty)$, $R \leq \frac{\rho^6}{16M^2}$, as well as $\varepsilon\in (0,1)$ let 
    Assumption~\ref{assum:Two} and Assumption~\ref{assum:coerc} be satisfied. Then 
    \begin{equation*}
         d_{\mathrm{H}} ( \lap ,\post ) =  \left(\int_{\R^d} \Big[\sqrt{\post(x)}-\sqrt{\lap(x)}  \Big]^2 \mathrm{d}x\right)^{1/2}  \in   \mathcal{O}\left(n^{-\frac{1}{2}}\right).
    \end{equation*}
\end{theorem}

The idea of the proof of this theorem follows the classical approach of `splitting the integral' appearing in the definition of the Hellinger distance into several parts and bound those suitably. 
The corresponding parts are determined by decomposing $\R^d$ into `high-density regions' around the local minima of $\logpostd$ as well as the remaining `tail region'. Under our assumptions the latter one is negligible.
The stated result exhibits the same rate of convergence for the Hellinger distance as \cite[Theorem~2]{art:SchillingBjornPhilipConvergLaplaceApprox} in a multi-modal setting. It is asymptotic in the sense that it relies on Lebesgue's dominated convergence theorem. Furthermore, we see in the proof that the dependence on $d$ of the bound on the Hellinger distance is exponential in this setting.

\subsection{Non-asymptotic result for $q \in \left[\frac{1}{3}, \frac{1}{2}\right)$}
From a theoretical and applied perspective, for example in the Bayesian inverse problem framework, it is particularly interesting, on the one hand, to have a fully explicit non-asymptotic estimate, and on the other hand, to achieve a `nice' dependence on $d$. To this end, 
we exploit stronger assumptions that enable us to indeed `shrink' the neighborhood $\Vnq$ in $n$ with rate $q \in \left[\frac{1}{3}, \frac{1}{2}\right)$. However, also the `tail area' needs to be controlled more explicitly, which motivates the following condition.
 
 \begin{assumption}\label{assum:Twobis}
 For $R,q,\delta>0$ assume that for all $n \in \N$, and $x \in \overline{V}_n^{R,q}$ holds 
\begin{equation*}
    \mathrm{e}^{-n\logpostd(x)} \le \sum_{k=1}^K \pi_n^{(k)} \exp \big( -n \delta \| x - x_n^{(k)} \|^2_{H_n^{(k)}} \big)  . 
\end{equation*}
\end{assumption}

\begin{rem}
Intuitively, \cref{assum:Twobis} is more restrictive than \cref{assum:coerc}, since it essentially imposes at least `Gaussian decay' at infinity, whereas \cref{assum:coerc} also allows for heavier `tail behaviors'. For instance, for the  $\logpostd\colon \R\to [0,\infty)$ with $\logpostd(x)  = \sqrt{x^2 + 1} -1$ there is an $\varepsilon\in (0,1)$, such that \cref{assum:coerc} is satisfied, but there are no numbers $R,q,\delta \in (0,\infty)$ such that \cref{assum:Twobis} is true.
\end{rem}

Note that the combination of \cref{assum:Two} and \cref{assum:Twobis} is inspired from \cite[Assumption~3.3]{art:HelinKretschmannNonAs} and extended to 
the multi-modal posterior distribution framework. It particularly allows us to obtain an explicit non-asymptotic result, in contrast to the asymptotic statement of \cref{theo:MainpropHelling2}.
We provide this result now and present its proof in Appendix~\ref{sec:proof}.

\begin{theorem} \label{theo:MainpropHelling} 
For $q \in \left[\frac{1}{3}, \frac{1}{2}\right)$ and $\rho,M,\delta \in (0,\infty)$ as well as $R = \rho^2/ M^{\frac{2}{3}}$ let \cref{assum:Two} and \cref{assum:Twobis} be satisfied.
Then, for each $n \in \N$ we have
    \begin{equation*}
 d_{\mathrm{H}}( \lap ,\post ) 
\leq \frac{13 \, d^{\frac{3}{2}}M}{n^{\frac{1}{2}} \rho^3} 
+ 2^{\frac{d}{2}+2}  \exp\Big(-\frac{\rho^2 \min\{1,2\delta\}}{8M^{\frac{2}{3}}} n^{(1-2q)}\Big) [(K+1) + \max\{1,(2\delta)^{-\frac{d}{2}}\}].
    \end{equation*}
\end{theorem}
A consequence of the previous theorem is the following illustrative result, where we provide the (technical) proof in Appendix~\ref{subsec: proof_coro}.
\begin{cor}
\label{cor: accesable_statement}
Under the assumptions of \cref{theo:MainpropHelling} for $\delta\geq 1/2$ the bound simplifies to
\begin{equation}
\label{eq: conseq_of_thm}
 d_{\mathrm{H}}( \lap ,\post ) 
\leq \frac{13 \, d^{\frac{3}{2}}M}{n^{\frac{1}{2}} \rho^3} 
+ 2^{\frac{d}{2}+2} (K+2) \exp\Big(-\frac{\rho^2 n^{(1-2q)} }{8M^{\frac{2}{3}}}\Big),
\end{equation}
such that $ d_{\mathrm{H}}( \lap ,\post ) \leq \frac{14 \, d^{3/2}M}{n^{1/2} \rho^3} $ for
\begin{equation}
\label{eq: n_large}
 n \geq \max\bigg\{ \Big [\frac{16 M^{2/3}}{\rho^2(1-2q)} \log \big( \frac{16 M^{2/3}}{\rho^2(1-2q)} \big)\Big]^{1-2q},\Big[\frac{16 M^{2/3}}{\rho^2} \log\big( \frac{2^{\frac{d}{2}+2}(K+2)\rho^3}{M} \big) \Big]^{\frac{1}{1-2q}} 
 \bigg\}.
\end{equation}
\end{cor}
Let us interpret the latter bound regarding the dependence on the dimension $d$ and $n$.
Think of $q \in \left[\frac{1}{3}, \frac{1}{2}\right)$ and 
$\rho,M,\delta \in [\frac{1}{2},\infty)$ being fixed. Then, from a rough estimation of \eqref{eq: n_large} we see that there is a constant $c>0$ (depending only on $q,\rho,M,\delta$) such that for any $\varepsilon\in(0,1)$ for $n\geq c(\,d^{3}\varepsilon^{-2} + d^{1/(1-2q)})$ we have $d_{\mathrm{H}}( \lap ,\post ) \leq \varepsilon$, which delivers a polynomial dependence on the dimension.

Another consequence of the previous theorem is an explicit bound on $\chi^2(\post \parallel \lap) $ by leveraging \eqref{eq:chi2_hellinger} and an auxiliary estimate (also used in the proof of Theorem~\ref{theo:MainpropHelling2}, see Lemma~\ref{lem:Mainlem},) to control the appearing ratio of normalizing constants.
It is proven in \Cref{sec:proofChiSq}.

\begin{cor}\label{theo:PropChiSq}
Under the assumptions of \cref{theo:MainpropHelling} for $\delta\geq 1/2$ and 
\begin{equation}
\label{eq: restr_chi_squ}
n \geq \max\bigg\{ \frac{8^2\,\mathrm{e^2}d^3 M^2}{\rho^6}, \Big[\frac{4M^{2/3}}{\rho^2} \log(2^{2+\frac{d}{2}}(K+2)) \Big]^{\frac{1}{1-2q}}\bigg\}
\end{equation}
we have $\frac{\post(x)}{\lap(x)} \le 2\mathrm{e}$ as well as
    \begin{equation*}
         \chi^2(\post \parallel \lap) \underset{\eqref{eq:chi2_hellinger}}{\leq} 12 d_{\mathrm{H}}( \lap ,\post )^2
         \underset{\eqref{eq: conseq_of_thm}}{\leq} 4056\,\frac{d^3M^2}{n \rho^6} 
+ 2^{d+9} (K+2)^2 \exp\Big(-\frac{\rho^2 n^{(1-2q)} }{4M^{\frac{2}{3}}}\Big).
    \end{equation*}
\end{cor}


We finish this section with a remark on the idea of proof of Theorem~\ref{theo:MainpropHelling2}.
\begin{rem}\label{rem:ideaproof}
In contrast to Theorem~\ref{theo:MainpropHelling2} in Theorem~\ref{theo:MainpropHelling} we use that the `tail' of the distribution is comparable to that of a normal distribution, cf. \cref{assum:Twobis}. 
With that and requiring for increasing $n$ a decrease of the `radius or size' of $\Vnq$  allows us to achieve the `nice' $d$ dependence in the bound above. Intuitively, in this situation the target measure becomes more concentrated as $n$ increases. Namely,
if the `radius or size' decreases as $n^{-q}$, then the regions of high probability of $\post$ remain contained in the $\Vnq$'s. 
\end{rem}

\section{Laplace-based Monte Carlo for expectations}

Now we are interested in the computation of the expectation of a quantity of interest $f\colon \R^d\to\R$ w.r.t. $\post$. For this we assume that $f$ is $\post$-integrable and all appearing random variables are defined on a sufficiently rich probability space $(\Omega,\mathcal{A},\mathbb{P})$. Moreover, let $\nu$ be a probability distribution on $\mathbb{R}^d$ which dominates $\post$, i.e., the density $\impw:=\frac{\mathrm{d} \post}{\mathrm{d} \nu}$ exists. Then, we can write
\begin{equation*}
    \mathbb{E}_{\post}[f] :=
    \int_{\R^d} f \,\mathrm{d} \post
    = \int_{\R^d} f(x)\ \impw(x)\,\nu(\mathrm{d}x).
\end{equation*}
Unfortunately, it is quite common that $\impw$ can only be evaluated up to a normalizing constant, i.e., it is only `partially known'. 
To take this into account, let $\unimpw\colon \R^d \to [0,\infty)$ be a function that returns values of 
$\impw$
up to a common multiplicative factor, such that
\begin{equation}
\label{eq:integr_problem}
\mathbb{E}_{\post}[f] = \frac{\int_{\R^d} f(x) \unimpw(x)\nu(\mathrm{d}x)}{\int_{\R^d} \unimpw(x)\nu(\mathrm{d}x)}.
\end{equation}

\subsection{Laplace-based importance sampling}\label{sec:LAIS}

For $n\in\mathbb{N}$ and an i.i.d. sequence $(X_{i,n})_{i\in\mathbb{N}}$ of
$\nu$-distributed random variables, define 
\emph{self-normalized importance sampling (SNIS)} w.r.t. $\nu$ as
\begin{equation*}
    \mathrm{IS}^{N}_{\post, \nu} (f):=\frac{ \sum_{i=1}^N \unimpw(X_{i,n}) \  f(X_{i,n})}{\sum_{j=1}^N  \unimpw(X_{j,n})}.
\end{equation*}
It is well-known that the strong law of large numbers yields consistency of SNIS, see \cite[Theorem 9.2]{mcbook}, i.e., we have for any $n\in\N$ that
\begin{equation*}
      \mathrm{IS}^{N}_{\post, \nu} (f) \underset{N \to \infty}{\overset{\text{a.s}}{\longrightarrow}} \mathbb{E}_{\post}[f] .
\end{equation*}
Moreover, by the central limit theorem 
and Slutzky's lemma one can prove the asymptotic normality of the SNIS estimator, i.e., we have
\begin{equation*}
    \sqrt{N}\left(\mathrm{IS}^{N}_{\post, \nu} (f) -  \mathbb{E}_{\post}[f] \right)  \underset{N \to +\infty}{\overset{\text{d}}{\longrightarrow}} \mathcal{N} (0 , \sigma_n^2(f) )
\end{equation*}
where $\overset{\text{d}}{\longrightarrow}$ denotes convergence in distribution, given the \emph{asymptotic variance} of SNIS is finite, that is,
\begin{equation}\label{eq:SNIS_AsympVar}
      \sigma_n^2(f) 
      :=
      \E_{\nu} \left[\impw^2 \left(f - \E_{\post}[f]\right)^2 \right]
      < \infty.
\end{equation}
All of these statements also apply if $\nu$ depends on $n$, i.e., if we consider a sequence of importance distributions $\nu_n$ to integrate w.r.t.~$\post$.

The asymptotic variance serves as a measure of the efficiency of SNIS. 
Although it depends on $f$, for bounded $f$ it is controlled by the second moment of the importance weights. Specifically, if $\|f\|_\infty < \infty$, then $\sigma_n^2(f)\leq 4 \|f\|_\infty^2\ \E_{\nu} \left[\impw^2\right]$, where $\E_{\nu} \left[\impw^2\right]$ coincides with $\chi^2(\post \parallel \nu) + 1$, see \cite[Section 2.3]{Agapiou2017Importance} .

This motivates another common measure of efficiency for SNIS: the \emph{effective sample size (ESS)}. 
It is defined as
\begin{equation}\label{eq:ESS}
    \ESS_n(N) 
    := 
    \frac{\left(\sum_{j=1}^N \unimpw(X_{j,n})\right)^2}{\sum_{i=1}^N \unimpw(X_{i,n})^2}.
\end{equation}
Thus, $\ESS_n(N)$ itself is a random variable that obviously satisfies $\ESS_n(N)\leq N$ almost surely.
As noticed in \cite{Agapiou2017Importance}, the expression $\frac{N}{\ESS_n(N)} - 1$ may serve as a proxy of $\chi^2(\post\parallel \nu)$. In fact, for $\chi^2(\post \parallel \nu) < \infty$ we have 
\begin{equation*}
      \frac{\ESS_n(N)}{N} \
      \underset{N \to \infty}{\overset{\text{a.s}}{\longrightarrow}} \
       \frac1{\chi^2(\post \parallel \nu) + 1}
\end{equation*}
and therefore for bounded $f$ and $\chi^2(\post \parallel \nu) < \infty$ we can control 
\[
    \sigma_n^2(f)
    \leq 4 \|f\|_\infty^2 \lim_{N\to\infty} \left(\frac{\ESS_n(N)}{N}\right)^{-1}. 
\]
Hence, $\frac1{\chi^2(\post \parallel \nu) + 1}$ can be viewed as \emph{asymptotic efficiency of SNIS} as the limit of the effective sample size and we set
\[
\mathrm{eff}^\mathrm{IS}_{\post}(\nu)
:=\frac1{\chi^2(\post \parallel \nu) + 1} \in[0,1].
\]

\paragraph{Prior-based importance sampling.}
We first consider $\nu=\pi_0$, i.e., the case of SNIS w.r.t. the prior distribution. We show that the corresponding asymptotic efficiency deteriorates as $n\to\infty$. For simplicity we restrict the statement to the case of unimodal $\post$.
\begin{theorem} \label{th: prior}
    Let $\Phi\colon \R^d \to [0,\infty)$ be a function independent of $n$ and define $I_n(x) := \Phi(x)-\frac{1}{n}\log(\prior(x))$ for any $n\in\N$ and all $x\in\R^d$. Moreover, assume that
    \begin{enumerate}
    \item 
    there is a unique global  
    minimizer $x_* \in \R^d$ of $\Phi$ with positive definite Hessian $H_* \coloneq \nabla^2 \Phi(x_{*})$; as well as
    \item
    for every $r > 0$ there holds
    $\inf_{x\in B_r^c(x_*)} \Phi(x) > \Phi(x_*),$
    where $B_r(x_*):=\{x\in\mathbb{R}^d:\|x-x_*\|\leq r\}$ and $ B_r^c(x_*):=\mathbb{R}^d \setminus B_r(x_*)$.
    
    \end{enumerate}
    Then, we have ${\mathrm{eff}^\mathrm{IS}_{\post}(\pi_0)} \in \mathcal{O}\left(n^{-d/2}\right).$
\end{theorem}
\begin{proof}
    We apply the classical Laplace method for the asymptotics of integrals \cite[Section~IX.5]{wong2001asymptotic} twice 
    and obtain as $n\to\infty$ that
    \begin{align*}
    \int_{\mathbb{R}^d} \exp\left(-n \Phi(x)\right) \pi_0(x) {\rm d}x 
    & 
    = \mathrm{e}^{-n\Phi(x_*)} n^{-d/2}\left(c_0(\pi_0) + \mathcal{O}\left(n^{-1}\right)\right), \quad \text{and} \\
\int_{\mathbb{R}^d} \exp\left(-2n \Phi(x)\right) \pi_0(x) {\rm d}x 
    &
    = \mathrm{e}^{-2n\Phi(x_*)} (2n)^{-d/2}\left(c_0(\pi_0) + \mathcal{O}\left(n^{-1}\right)\right),
    \end{align*}
    where $c_0(\pi_0):=\sqrt{\det(2\pi H_*^{-1})}\pi_0(x_*)>0$ by the assumptions.
    Consequently, we obtain by $Z_n= \int_{\mathbb{R}^d} \exp\left(-n \Phi(x)\right) \pi_0(x) {\rm d}x $ and by the previous representations that 
    \begin{align*}
        \mathrm{eff}^\mathrm{IS}_{\post}(\pi_0) 
        &
        = \frac{1}{\E_{\pi_0}(\impw^2)} 
        = \frac{\Big( \int_{\mathbb{R}^d} \exp\left(-n \Phi(x)\right) \pi_0(x) {\rm d}x   \Big)^2}{\int_{\R^d} \exp\left(- 2n \Phi(x)\right) \pi_0(x) {\rm d}x}
        = \frac{n^{-d}\left(c_0(\pi_0) + \mathcal{O}\left(n^{-1}\right)\right)^2}{(2n)^{-d/2}\left(c_0(\pi_0) + \mathcal{O}\left(n^{-1}\right)\right)},
    \end{align*}
    which implies the claim.
\end{proof}
This theorem shows that the prior-based importance sampling suffers from a curse of dimensionality and also becomes less efficient as the target distribution $\pi_n$ becomes more concentrated.

\paragraph{Laplace-based importance sampling.}
Next we consider SNIS based on the multi-modal version of the LA, i.e., $\nu=\lap$. In this setting the corresponding weights are given by
\begin{equation}
    \unimpw(x) 
    := 
    \frac{\exp\left(-n I_n\left(x\right)\right)}{\sum_{k=1}^K \pi_n^{(k)}\exp\big( - \frac{n}{2} \|x-x_n^{(k)}\|^2_{H_n^{(k)}}\big)},
    \qquad
    x\in\R^d.
    \label{weight_LM}
\end{equation}


We start with statements about the asymptotic efficiency $\mathrm{eff}^\mathrm{IS}_{\post}(\lap)$ of Laplace-based importance sampling, which is
an immediate consequence of Corollary~\ref{theo:PropChiSq}.
\begin{cor}\label{theo:chi2}
    For $q \in \left[\frac{1}{3}, \frac{1}{2}\right)$, $\delta\geq1/2$ and $\rho,M \in (0,\infty)$ as well as $R = \rho^2/ M^{\frac{2}{3}}$ let \cref{assum:Two} and \cref{assum:Twobis} be satisfied.
    Then, for Laplace-based importance sampling we have 
    \[
        1 - \mathrm{eff}^\mathrm{IS}_{\post}(\lap)
        \in \mathcal O(n^{-1}).
    \]
\end{cor}
However, particularly, \cref{assum:Twobis} is not necessarily satisfied in Bayesian inference with $\delta \geq \frac 12$.
Therefore, we next extend our analysis of SNIS based on the multimodal version of the Laplace approximation and study the behavior of $\ESS_n(N)$ under \cref{assum:Two}. For this, we formulate the following additional requirement.


\begin{assumption} \label{assum:minimizer-limit}
Let $\Phi\in \mathcal C^2(\R^d, [0,\infty))$ be a function that is independent of $n$ and define
    $I_n(x) := \Phi(x) - \frac{1}{n} \log(\prior(x))$ for any $n\in\N$ and all $x\in\R^d$. Assume that
    $\Phi$ has $K\in\mathbb{N}$ global minimizers $x_{*}^{(1)},\dots,x_{*}^{(K)} \in \R^d$ satisfying
	$\Phi (x_{*}^{(k)}) = 0$ with positive definite $\nabla^2 \Phi (x_{*}^{(k)}) $ and
    \[
        \lim_{n\to\infty} x_{n}^{(k)} = x_{*}^{(k)},
        \qquad k = 1,\ldots,K.
    \]
\end{assumption}

\begin{rem}
       The requirement that $\Phi (x_{*}^{(k)})=0$ for any $k\in\{1,\dots,K\}$ is not restrictive. Namely, if $\min_{z \in \R^d} \Phi (z)> 0$,  then we can simply set
    \[
        \widehat{\Phi} \left(x\right) \coloneq \Phi \left(x\right) - \min_{z \in \R^d} \Phi \left(z\right), \qquad \widehat{I}_n(x) \coloneq \widehat{\Phi}(x) - \frac{1}{n} \log(\prior(x)),
    \]
    for which we obtain $\min_{z \in \R^d} \widehat{\Phi} \left(z\right) = 0$ as well as
    \[
    \post(\mathrm{d}x) = \frac{1}{\widehat{Z}_n}\mathrm{e}^{-n \widehat{\logl}_n(x)} \pi_0(\mathrm{d}x), 
    \quad\text{with}\quad 
    \widehat{Z}_n \coloneq  \int_{\R^d} \mathrm{e}^{-n \widehat{\logl}_n(x)}\prior(\mathrm{d}x).
    \]
    However, the assumption that $\lim_{n\to\infty} x_{n}^{(k)} = x_{*}^{(k)}$ requires that all local minima of $I_n$ converge to global minima of $\Phi$ which does, for instance, exclude the case of vanishing local modes of $\post$.
\end{rem}

\begin{rem}
    It is worth to emphasize here that the conditions of \cref{assum:minimizer-limit} (and also the ones of \cref{th: prior}) are in line with our formulated standing assumption from Section~\ref{sec:setting}. These just provide a refinement of requirements that can be further exploited and allow new conclusions. 
\end{rem}


In our (asymptotic) analysis of the effective sample size of SNIS we exploit the fact that for increasing `concentration' (i.e. increasing $n$) iid samples w.r.t. $\lap$ increasingly `concentrate' around modes. We have the following result which is proven in Section \ref{sec:GnN_proof}.


\begin{prop}\label{prop: GnN}
Let $R>0$, $N\in\N$, $q\in[0,1/2)$ and $X_{1,n},\dots,X_{N,n}$ be iid $\lap$-distributed random variables as well as
 \begin{equation}\label{eq:GnN}
        G_{n,N} 
        \coloneqq 
        \bigcap_{i=1}^N 
        \Big\{ 
        \min_{k=1,\ldots,K} 
        \big\| X_{i,n} - x_n^{(k)} \big\|^2_{H_n^{(k)}} 
        \leq 
        Rn^{-2q} 
        \Big\}.
    \end{equation}
    Then, $\lim_{n \to \infty} \mathbb{P}(G_{n,N}) = 1.$
\end{prop}

The proof can be found in Section \ref{sec:GnN_proof}. 
As a consequence of these former result we are able to add the following tool which is proven in Section~\ref{app:lem W*}.
\begin{lem}\label{lem: W*}
For $q \in \left(\frac{1}{3}, \frac{1}{2} \right)$ and $R,\rho,M\in (0,\infty)$
let \cref{assum:Two} as well as \cref{assum:minimizer-limit} be satisfied.
For $N\in\mathbb{N}$ let $X_{1,n},\dots,X_{N,n}$ be iid $\lap$-distributed random variables and
$G_{n,N}$ be defined as in \eqref{eq:GnN}.
    Then, there exists a constant $c>0$ and an $n_0 \in \mathbb{N}$, such that for all $n\in\mathbb{N}$ with $n \ge n_0$ for $\ESS_n(N)$ (based on $X_{1,n},\dots,X_{N,n}$) restricted\footnote{
    For a real-valued random variable $Z$ (defined on $(\Omega,\mathcal{A},\mathbb{P})$) a statement (for example an inequality), restricted to $G\in\mathcal{A}$, holds if it is true for any $\omega\in G$.} 
    to $G_{n,N}$ holds
    \begin{equation}
    \label{eq: ESS-N}
   \left| \ESS_n(N) - N \right| \leq c N n^{2-6q}.
   \end{equation}
\end{lem}

Combining the former two results leads to the following asymptotic statement regarding the effective sample size.

\begin{theorem}\label{theo: ESS}
For $q \in \left(\frac{1}{3}, \frac{1}{2} \right)$ and $R,\rho,M\in (0,\infty)$
let \cref{assum:Two} as well as \cref{assum:minimizer-limit} be satisfied.
For $N\in\mathbb{N}$ let $X_{1,n},\dots,X_{N,n}$ be iid $\lap$-distributed random variables. Then, for $\ESS_n(N)$ (based on $X_{1,n},\dots,X_{N,n}$) and any $
\gamma \in (0,6q-2)$ we have
\[
n^{\gamma}\left| \ESS_n(N) - N \right| 
\xrightarrow[n\to\infty]{\mathbb P}
0,
\]
where $\xrightarrow[n\to\infty]{\mathbb P}$ denotes convergence in probability.
\end{theorem}
\begin{proof} Let $\varepsilon >0$ be arbitrary. For $c>0$ and $n_0$ from Lemma~\ref{lem: W*}, for $n\in\mathbb{N}$ with $n\geq n_0$ and $n\geq (cN/\varepsilon)^{1/(6q-2-\gamma)}$, restricted to $G_{n,N}$ (defined as in \eqref{eq:GnN}), we have
\[
 n^{\gamma}\left| \ESS_n(N) - N \right| 
 \underset{\eqref{eq: ESS-N}}{\leq}  cN n^{\gamma+2-6q} \leq \varepsilon.
\]
Therefore
\begin{align*}
        \lim_{n\to\infty} \mathbb{P} \left( n^\gamma \left| \ESS_n(N) - N \right| \leq \varepsilon
        \right)
        &\geq \lim_{n\to\infty} \mathbb{P} \left( n^{\gamma} | \ESS_n(N) - N | \leq \varepsilon \text{ and } G_{n,N}\right)
        = \lim_{n\to\infty}  \mathbb{P} ( G_{n,N} ).
\end{align*}
By \cref{prop: GnN} holds $ \lim_{n\to\infty} \mathbb{P} \left( G_{n,N} \right) = 1$, which finishes the proof.
\end{proof}


\subsection{Laplace-based independent Metropolis-Hastings algorithm}\label{sec:LAIMH}
We consider an alternative to SNIS for the approximation of $\E_{\post}[f]$ that is based on Markov chain Monte Carlo, more specifically on the \emph{Independent Metropolis-Hastings} (IMH) algorithm. The idea is to construct a Markov chain that can be used for approximate sampling of $\post$. Recall that $\nu$ is a probability measure on $\R^d$, 
such that the density $\impw:=\frac{\mathrm{d} \post}{\mathrm{d} \nu}$ exists and is assumed to 
coincide up to a multiplicative constant with $\unimpw$. We consider the IMH algorithm based on proposal distribution $\nu$, which provides a construction for a $\post$-reversible Markov chain $(X_{i,n})_{i\in\N}$. 
Namely, for $i\in\N$ given $X_{i,n} = x\in \R^d$ the next instance $X_{i+1,n}$ is specified by the following transition mechanism: 
First, draw $Y\sim\nu$ and call the result $y$, afterwards with \emph{acceptance probability}
\[
    \alpha_n(x,y) 
    = 
    1 \wedge \frac{\unimpw(y)}{\unimpw(x)},
\]
accept the proposal $y$, i.e., set $X_{i+1,n} = y$ and otherwise reject it, i.e., set $X_{i+1,n} = x$.
The corresponding transition kernel, that determines the conditional distribution of $X_{i+1,n}$ given $X_{i,n}=x$, can be represented as 
\begin{align*}
    \mathrm{IM}_n(x,\mathrm dy) := \alpha_n(x,y) \nu(\mathrm dy) + r_n(x) \delta_x(\mathrm dy),
\end{align*}
where $r_n(x) := \int_{\R^d} \left( 1-\alpha_n(x,y) \right) \nu(\mathrm dy)$ denotes the rejection probability and $\delta_x$ the Dirac measure at $x$.
A sequence of $\R^d$-valued random variables $(X_{i,n})_{i\in\N}$ that can be realized iteratively by the transition mechanism is called a Markov chain with transition kernel $\mathrm{IM}_n$.
In our setting this Markov chain is ergodic\footnote{By the absolute continuity of $\nu$ w.r.t. $\post$, it is easily seen that the IMH algorithm satisfies for any Borel measurable $A\subseteq \R^d$ with $\post(A)>0$ that $\mathrm{IM}_n(x,A)>0$ for any $x\in \mathbb{R}^d$. This readily implies $\post$-irreducibility, aperiodicity and by \cite[Corollary~2]{tierney1994markov} Harris recurrence, such that by \cite[Theorem~1]{tierney1994markov} the ergodicity statement follows.}, i.e., we have
\[
    \lim_{t\to\infty} d_{\mathrm{TV}} \left( \mathrm{IM}_n^t(x,\cdot),\post \right)
    =
    0
    \qquad
    \forall x\in\R^d.
\]
Moreover, for any $\post$-integrable $f\colon \R^d \to \R$ a strong law of large numbers applies, see \cite[Theorem~3]{tierney1994markov}, such that
\begin{equation*}
      \mathrm{IMH}^N_{\post,\nu}(f)
      \underset{N \to \infty}{\overset{\text{a.s}}{\longrightarrow}} \mathbb{E}_{\post}[f],
      \qquad
         \mathrm{IMH}^N_{\post,\nu}(f) := 
   \frac{1}{N} \sum_{i=1}^N f(X_{i,n}).
\end{equation*}
If the function $\impw$ is bounded, i.e., 
\begin{equation}\label{eq:w_bounded}
 \| \impw \|_{\infty} \coloneqq \sup_{x\in\mathbb{R}^d} \impw(x)
 < \infty,
\end{equation}
then the IMH Markov chain is uniformly ergodic (cf.~\cite{MengersenTweedie1996}) and, by \cite[Theorem~2]{wang2020exact}, we have
\begin{equation}
    \label{eq:worst_case_tv}
        \sup_{x\in\R^d}
    d_{\mathrm{TV}} \left( \mathrm{IM}_n^t(x,\cdot),\post \right)
    =
    \left(1 - \frac{1}{\| \impw \|_{\infty}}\right)^t,
    \qquad
    t\in\N.
\end{equation}
In this case, there also holds a central limit theorem (cf. \cite[Corollary 2.1]{RoRo1997}) for $\mathrm{IMH}^N_{\post,\nu}(f)$, i.e., 
\begin{equation*}
    \sqrt{N}\left(\mathrm{IMH}^{N}_{\post, \nu} (f) -  \mathbb{E}_{\post}[f] \right)  \underset{N \to \infty}{\overset{\text{d}}{\longrightarrow}} \mathcal{N} (0 , \widetilde \sigma_n^2(f) ),
    \qquad
    f \in L^2(\post),
\end{equation*}
where the \emph{asymptotic variance} is given by
\begin{equation}\label{eq:IMH_AsympVar}
      \widetilde \sigma_n^2(f) 
      :=
      \Var_{\post}(f)
      \Big(1 + 2\sum_{i=1}^\infty \mathrm{Corr}\big(f(\widetilde X_{1,n}), f(\widetilde X_{i,n}) \big)  \Big),
\end{equation} 
for $(\widetilde X_{i,n})_{i\in\N}$ denoting an IMH Markov chain with initial distribution $\post$ and $\mathrm{Corr}(\cdot,\cdot)$ denoting the correlation between random variables. Besides that, the asymptotic variance also characterizes the asymptotic mean squared error, that is, for example by virtue of \cite[Theorem~3.34]{R12}, we have
\[
 \lim_{N\to\infty} N\cdot \E \big\vert  \mathrm{IMH}^{N}_{\post, \nu} (f) - \E_{\post}[f] \big\vert^2 = \widetilde \sigma_n^2(f).
\]
Therefore, we are interested in the $n$-dependence of this quantity. 
Note that the asymptotic variance is closely related to spectral properties of the Markov operator
$M_n\colon L^2(\post) \to L^2(\post)$, associated with the transition kernel $\mathrm{IM}_n$, given by
\[
    (M_n f)(x) := \int_{\mathbb{R}^d} f(y)\,\mathrm{IM}_n(x,{\rm d}y), \quad f\in L^2(\post).
\]
Define $L^2_{0}(\post)\subseteq L^2(\post)$ as the space of functions $f\colon \R^d \to \R$ satisfying $\E_{\post}[f]=0$.
Then, the $L^2(\post)$-spectral gap of $M_n$ is given as
\[
    \mathrm{gap}_{\post}(M_n) := 1 - \|M_n\|_{\post},
    \qquad
    \|M_n\|_{\post}
    := \sup_{\substack{f \in L^2_{0}(\post), f \neq 0}}
    \frac{ \|M_n f\|_{L^2(\post)} }{ \|f\|_{L^2(\post)} }.
\]
Clearly, $\mathrm{gap}_{\post}(M_n)\in [0,1]$ and it is well known, see for example \cite[p.~6]{rudolf2018generalization}, that
\[
    \widetilde \sigma_n^2(f) \leq \frac{2 \Vert f\Vert_{2,\post}}{\mathrm{gap}_{\post}(M_n)},
    \qquad \forall f \in L^2(\post).
\]
Therefore, non-deteriorate $n$-dependence of $\mathrm{gap}_{\post}(M_n)$ delivers insights into non-deteriorate $n$-dependence of the asymptotic variance. Independent of that, the spectral gap itself is an important quantity in the investigation of stabilizing behavior of Markov chains, see for example \cite{R12} or the recent discussions in \cite{hofstadler2026solving}. For the spectral gap of the IMH Markov operator we obtain the following lower bound, which allows a similar conclusion as \eqref{eq:worst_case_tv}, cf. \cite{wang2020exact}.
\begin{prop} \label{prop:IMH_SpectralGap}
    If $\| \impw \|_{\infty}<\infty$, then
    \begin{equation}
        \label{eq:IMH_SpectralGap}
        \mathrm{gap}_{\post}(M_n) \geq \frac{1}{\| \impw \|_\infty}.
    \end{equation}
\end{prop}
In particular, by \cite[Theorem 2.1]{MengersenTweedie1996}, 
if $\| \impw \|_{\infty} = \infty$, then an IMH Markov chains fails to be geometrically ergodic. Consequently, by virtue of \cite{gallegos2024equivalences}, in that case $\mathrm{gap}_{\post}(M_n)=0$, such that w.r.t. the spectral gap of IMH there is (as visible in \eqref{eq:worst_case_tv} for the total variation convergence) a dychotomy of its existence, determined by whether  
$\| \impw \|_{\infty}$ is finite or not.
The result of the proposition and the former conclusion is actually well known, cf. \cite{atchade2007geometric,gaasemyr2006spectrum,MengersenTweedie1996,wang2020exact}, but in our setting not really explicitly  stated, such that we provide a proof in \Cref{app:IMH_SpectralGap}.

Besides the asymptotic variance and/or spectral gap we provide seemingly less prominent, easier to handle quantities, that, however, have less direct consequences for stabilization. 
One is the \emph{averaged acceptance probability} of IMH w.r.t. proposal distribution $\nu$, given as
\begin{equation}
\mathbb{E}_{\post\otimes\nu}[\alpha_n] := \int_{\R^d}\int_{\R^d} \alpha_n(x,y) \ \nu(\mathrm dy) \post(\mathrm d x) \in[0,1],
\label{eq:avg_ar}
\end{equation}
and the other one is the \emph{expected squared jump distance in a direction} $v\in\R^d$ given as
\begin{equation}
        \rho_n^v := \int_{\R^d}\int_{\R^d} \left|v^Tx-v^Ty\right|^2 \alpha_n(x,y)\nu({\rm d}y)\pi_n({\rm d}x).
\label{eq:rho_v_n}
\end{equation}
Note, that $\rho_n^v$ relative to the variance of the target $\post$ in the direction $v$ corresponds directly to the lag-one autocorrelation of an IMH Markov chain $(\widetilde X_{i,n})_{i\in\N}$, with proposal $\nu$ and initial distribution $\post$, in this direction, i.e., we have
\begin{equation}
    \mathrm{Corr}(v^\top \widetilde{X}_{i,n}, v^\top \widetilde{X}_{i+1,n})
    =
    1 - \frac{\rho_n^v}{2\Var_{X\sim\post}(v^TX)}.
\label{eq:Corr_rho_v_n}
\end{equation}
In the case of an averaged acceptance probability and/or an expected squared jump distance that is zero the convergence behavior of the IMH Markov chain is `bad'. Namely, on average no proposed state is accepted and/or one stays always at the current state and got stuck. Therefore, asymptotic in $n$ we desire a constantly large averaged squared jump distance as well as averaged acceptance probability strictly larger than zero. For statements regarding these two simple measures of efficiency the following tool
is useful, which is proven in Appendix~\ref{app:EAA}.
\begin{lem}\label{lem:EAA}
    For approximate sampling of $\post$
    with the IMH algorithm with proposal $\nu$ the average acceptance probability satisfies
    \[
    1 - 2 d_{\mathrm{TV}}(\post,\nu)
    \leq
    \mathbb{E}_{\post\otimes\nu}[\alpha_n]
    \leq
    1 - d_{\mathrm{TV}}(\post,\nu).
    \]
    Moreover, if $\nu$ and $\post$ possess 4th moments, then for any $v \in \R^d$, 
    with $(\widetilde X_{i,n})_{i\in\N}$ being an IMH Markov chain with proposal $\nu$ and initial distribution $\post$,
    we have 
    \[
        \left| \mathrm{Corr}(v^\top \widetilde{X}_{i,n}, v^\top \widetilde{X}_{i+1,n})\right|
        \leq
        \sqrt{\frac{2 \int_{\R^d}\int_{\R^d} |v^Tx-v^Ty|^4 
        \left[\nu(\mathrm dy) + \post(\mathrm d y)\right]
        \post(\mathrm d x)}{\Var_{X\sim \pi_n}(v^\top X)^2}}\ 
        d_{\mathrm H} \left(\pi_n,\nu \right).
    \]
\end{lem}



\paragraph{Laplace-based IMH.} We now consider choosing the (multimodal) Laplace approximaton as proposal distribution $\nu= \lap$ in the IMH algorithm. As an immediate corollary of \Cref{theo:MainpropHelling2} or \Cref{theo:MainpropHelling}, together with \Cref{lem:EAA}, we obtain the following result.



\begin{theorem}\label{thm:LAIMH}
    Let the assumptions of \Cref{theo:MainpropHelling2} or \Cref{theo:MainpropHelling} be satisfied. 
    Then the average acceptance probability of Laplace-based IMH with proposal distribution $\lap$ satisfies
    $$
    1 - \mathbb{E}_{\post\otimes\lap}[\alpha_n] \in \mathcal{O}(n^{-1/2}).
    $$
    Moreover, if in addition 
    $\pi_0$ possesses a fourth moment and $\post$ concentrates around finitely many fixed minimizers $x_*^{(1)},\cdots,x_*^{(K)}$ with $K \geq 2$,
    then for any $v \in \mathbb{R}^d$ with $\lim_{n\to\infty} \Var_{X\sim \pi_n}(v^\top X) > 0$ we have for an IMH Markov chain $(\widetilde{X}_{i,n})_{i\in\N}$ with proposal distribution $\lap$ and initial distribution $\post$ that
\[
    \left|\mathrm{Corr}(v^\top \widetilde{X}_{i,n}, v^\top \widetilde{X}_{i+1,n}) \right|
    \in 
    \mathcal{O}(n^{-1/2})
\]
\end{theorem}
\begin{proof}
    The first statement is an immediate consequence of \Cref{lem:EAA} and equation \eqref{eq:dTV_dH}, such that
    $$1 - \mathbb{E}_{\post\otimes\lap}[\alpha_n]
    \leq
    2 d_{\mathrm{TV}}(\post,\lap)
    \leq 2 d_{\mathrm{H}}(\post,\lap),$$ which, combined with the statements of \Cref{theo:MainpropHelling2} or \Cref{theo:MainpropHelling}, respectively, verify the claim.
    The second statement follows analogously by \Cref{theo:MainpropHelling2} or \Cref{theo:MainpropHelling}, respectively, and noting that the finite fourth moment of $\pi_0$ guarantees that \[\int_{\R^d}\int_{\R^d} |v^Tx-v^Ty|^4 
        \left[\lap(\mathrm dy) + \post(\mathrm d y)\right]
        \post(\mathrm d x) 
        < \infty. \]
    Together with the posterior concentration assumption and either \cref{assum:coerc} or \cref{assum:Twobis}, this further implies that \[\int_{\R^d}\int_{\R^d} |v^Tx-v^Ty|^4 
        \left[\lap(\mathrm dy) + \post(\mathrm d y)\right]
        \post(\mathrm d x) 
        \in \mathcal{O}(1). \]
    Since $\lim_{n\to\infty} \Var_{X\sim\post}\left(v^T X\right) > 0$ and the fourth-order integral remains bounded as $n\to\infty$, the prefactor in the absolute lag-one correlation bound also remains bounded. Hence, the absolute lag-one correlation inherits the $\mathcal{O}(n^{-1/2}) $ upper rate of the Hellinger distance. \qedhere
\end{proof}

\begin{rem}
    The assumption in \Cref{thm:LAIMH} about the fourth moment of $\pi_0$ is essential to guarantee the finiteness of $\int_{\R^d}\int_{\R^d} |v^Tx-v^Ty|^4 
        \left[\lap(\mathrm dy) + \post(\mathrm d y)\right]
        \post(\mathrm d x)$. However, the statement can be generalized to cases where $\post$ is unimodel or $\lim_{n\to\infty} \Var_{X\sim \pi_n}(v^\top X) = 0$. The proof in these cases is more technical and exploits that the asymptotic behaviour of $\Var_{X\sim \pi_n}(v^\top X)$ is $\mathcal O(n^{-1/2})$ under suitable assumptions which then also yield that $\int_{\R^d}\int_{\R^d} |v^Tx-v^Ty|^4 
        \left[\lap(\mathrm dy) + \post(\mathrm d y)\right]
        \post(\mathrm d x)$ is of order $\mathcal O(n^{-1/2})$. We refer to \cite[Section 4.2]{rudolf2022robust} for more details in the unimodal case.
\end{rem}

Next, we establish an even stronger statement about the spectral gap of the Laplace-based IMH given the setting of \Cref{theo:PropChiSq}. Its proof relies also to some extent of the arguments delivered in \Cref{sec:proofChiSq}.

\begin{theorem}\label{theo:sp}
    Let the assumptions of \Cref{theo:PropChiSq} be satisfied with $\delta>1/2$.
    Then, for the spectral gap of the IMH algorithm with proposal distribution $\nu_n = \lap$, the associated Markov operator $M_n$ satisfies
    \[
        1 - \mathrm{gap}_{\pi_n}(M_n)\in \mathcal{O}\!\left(n^{1-3q}\right).
    \]
\end{theorem}
\begin{proof}
We provide an upper bound of $\| \impw \|_\infty$, from which the statement follows by Proposition~\ref{prop:IMH_SpectralGap}. 
Note that $\impw (x) = \frac{\widetilde Z_n^K}{Z_n}\,\unimpw(x) $. By \eqref{al:cite_for_spectral_gap}
 holds $\frac{\widetilde Z_n^K}{Z_n}\leq (1- \frac{\mathcal{K}_{n,1}}{2})^{-1}$, where $\mathcal{K}_{n,1}$ is defined in Lemma~\ref{lem:Mainlem}. Therefore we focus on bounding $\unimpw$ suitably. 
For $x \in \bigcup_{k=1}^K V_n^{k,R,q}$, by \eqref{eq:important_bnd_spectral_gap} we have
\[
    \unimpw \left(x\right)
    \leq
    \mathrm e^{n^{1-3q}}.
\]
In contrast to what we relied on in the proof of \Cref{theo:PropChiSq},
for $x\in\overline{V}_n^{R,q}$, by \cref{assum:Twobis}, we use
 \begin{align*}
     \unimpw (x) &\le     \sum_{k=1}^K \pi_n^{(k)}\frac{ \exp\big(-n \delta  \|x - x_n^{(k)} \|^2_{H_n^{(k)}}\big)}{\sum_{j=1}^K \pi_n^{(j)}  \exp\big(-\frac{n}{2} \|x - x_n^{(j)} \|^2_{H_n^{(j)}} \big)}
     \leq \sum_{k=1}^K  
     \mathrm{e}^{-n ( \delta-\frac{1}{2} )  \|x - x_n^{(k)} \|^2_{H_n^{(k)}}}
     \leq K \mathrm{e}^{-( \delta-\frac{1}{2} )R n^{1-2q}},
 \end{align*}
 such that 
\[
\sup_{x\in\mathbb{R}^d} \unimpw (x) 
\leq
\max\left\{ \mathrm e^{n^{1-3q}}, K \mathrm{e}^{-( \delta-\frac{1}{2} )R n^{1-2q}}\right\}.
\]
By $q \in [1/3,1/2)$, there exists a number $c>0$, such that for all sufficiently large $n\in\N$ holds 
$
    \| \unimpw \|_{\infty} \leq \mathrm e^{cn^{1-3q}}
$
and therefore $\| \impw \|_\infty \leq (1- \frac{\mathcal{K}_{n,1}}{2})^{-1} \mathrm e^{cn^{1-3q}}$.
By \eqref{eq:IMH_SpectralGap}, also for sufficiently large $n\in\N$, we obtain
\begin{align*}
1 - \mathrm{gap}_{\pi_n}(M_n) 
& \leq 1 - \frac{1}{\| \impw \|_{\infty}} 
\leq 1 -  \big(1- \frac{\mathcal{K}_{n,1}}{2}\big) \mathrm e^{-cn^{1-3q}} 
\leq 1 - \mathrm e^{-cn^{1-3q}} + \frac{\mathcal{K}_{n,1}}{2} \mathrm e^{-cn^{1-3q}}\\
&\leq 1 - \mathrm e^{-cn^{1-3q}} + \mathcal{K}_{n,1}
\leq cn^{1-3q} + \mathcal{K}_{n,1},
\end{align*}
where in the last equality we applied $1 - \mathrm e^{-x} \leq x$ which holds for all $x\in\R$.
By Lemma~\ref{lem:Mainlem} part \eqref{en: stat_b} we have $\mathcal{K}_{n,1} \in \mathcal{O}\!\left(n^{-1/2}\right)$. Consequently by $q \in [1/3,1/2)$ this yields $cn^{1-3q} + \mathcal{K}_{n,1} \in \mathcal{O}\!\left(n^{1-3q}\right)$ and the claim is proven.
\end{proof}


    \Cref{theo:sp} shows that the spectral gap of the Markov operator $M_n$ approaches one as $n\to\infty$, since $1-3q<0$ for $q\in(\frac{1}{3},\frac{1}{2})$. Hence, the convergence rate of the IMH algorithm improves asymptotically, indicating that the Laplace approximation $\lap$ becomes an increasingly effective proposal distribution for the target $\pi_n$. However, the condition $\delta>1/2$ is relatively restrictive and is essential for the present result. Without this assumption, the above argument no longer guarantees such asymptotically optimal behavior. In more general settings, IMH may exhibit subgeometric convergence. For instance, \cite{deligiannidis2024importance} establishes polynomial convergence rates in total variation under suitable moment assumptions. Such behavior, however, is not analyzed in the present work.





\section{Numerical Illustrations}\label{sec:numerics}

We now illustrate our theoretical findings on the Laplace-based importance sampling and Laplace-based independent Metropolis--Hastings algorithm at two examples. The first one is a rather academic test example from \cite{art:SchillingBjornPhilipConvergLaplaceApprox}. The second example is concerned with sampling the posterior in Bayesian neural networks where multimodal posteriors occur naturally by the symmetric architecture of the network. 

\subsection{Simple Example}\label{sec:simpleExample}
We consider a Gaussian prior $\pi_0 = \mathcal{N}(0,I_d)$ in $\mathbb{R}^d$ for various $d \in \{1,2,\ldots,7\}$ which is then conditioned on observations of
\[
    Y = \mathcal{G}(X) + \eta_n
\]
where $X\sim \pi_0$ and $\eta_n$ are stochastically independent and the forward map $\mathcal{G}(x) = (\mathcal{G}_1(x),\ldots,\mathcal{G}_d(x))^\top$ is given by
\begin{align*}
    &\mathcal{G}_1(x) = x_1^2 + x_1, \quad \mathcal{G}_2(x) = 2x_2 + x_1^2, \quad \mathcal{G}_3(x) = \exp(x_3)+x_2, \quad \mathcal{G}_4(x) = x_4+2x_3,\\
    &\mathcal{G}_5(x) = x_5^2 + x_4, \quad \mathcal{G}_6(x) = x_6+x_5, \quad \mathcal{G}_7(x) = x_7+2x_6
\end{align*}
for $x = (x_1,\ldots, x_d)^\top \in \R^d$. We consider two scenarios in which we want to sample the resulting posteriors $\post$ as given in \eqref{eq:definitionofmu}
\begin{enumerate}
    \item \emph{Small noise limit:} a single observation on increasing precision is given $y_n = \mathbf{2}_d + \eta_n$, $\mathbf{2}_d:=(2,\cdots,2)\in\mathbb{R}^d$, where $\eta_n \sim N(0, 0.1 \, n^{-1} \, I_d)$, and $\post$ follows \eqref{eq:definitionofmu} with
    \[
            \logl_n(x) := \frac{10}{2}\|y_n - \mathcal{G}(x)\|^2
    \]
    
    \item \emph{Large data limit:} an increasing number of observations of same precision are given, i.e., $y_i = \mathbf{2}_d + \eta_i$ for $i=1,\ldots,n$ with $\eta_i \sim N(0,0.1 I_d)$ being stochastically independent, hence $\post$ follows \eqref{eq:definitionofmu} with
    \[
            \logl_n(x) := \frac{10}{2n} \sum_{i=1}^n \|y_i - \mathcal{G}(x)\|^2.
    \]
\end{enumerate}
In each case the resulting posteriors are multimodal: for $d \in \{1,\ldots,4\}$ we have $K=2$ modes and for $d\in\{5,6,7\}$ we have $K=4$ modes. 
We apply LA-based importance sampling and LA-based IMH to sample from the resulting posterior $\post$ and investigate their statistical efficiency. 

\begin{figure}[H]
    \centering
    \begin{subfigure}[H]{0.4\textwidth}
        \includegraphics[width=\linewidth]{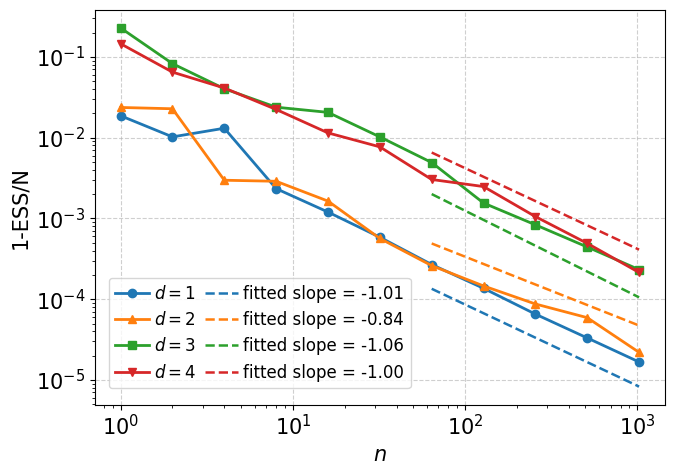}
        \caption{Small noise limit}
    \end{subfigure}
        \quad
    \begin{subfigure}[H]{0.4\textwidth}
        \includegraphics[width=\linewidth]{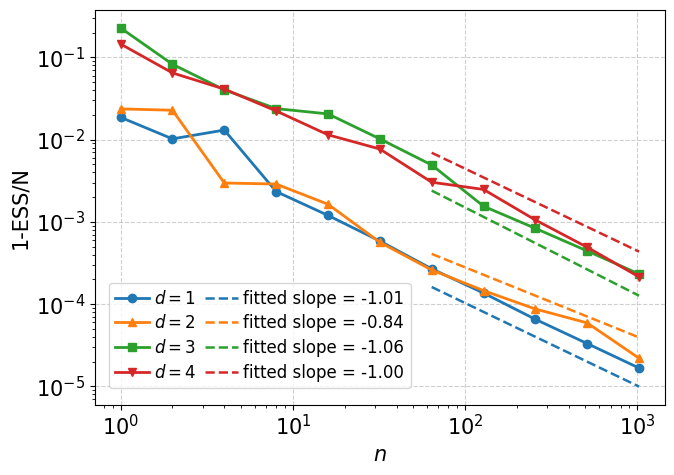}
        \caption{Large data limit}
    \end{subfigure}
    \caption{Performance of Laplace-based importance sampling for the simple example considered in Section \ref{sec:simpleExample}.}
    \label{IS}
\end{figure}

\paragraph{Laplace-based importance sampling:} We use $N=10^3$ samples for self-normalized importance sampling with $\nu_n = \lap$ in both scenarios. Figure \ref{IS} illustrate  the asymptotic behavior of the effective samples size for increasing $n$ and different dimensions $d$ for both scenarios. 
The numerical results show good agreement with our theoretical analysis: $1 - \frac{\ESS}{N}$ decays like $n^{-1}$ in the left panel, in accordance with Theorem \ref{theo: ESS}. Moreover, also in the large data limit we observe the same asymptotic behavior as in the small noise limit case, i.e., also here $1 - \frac{\ESS}{N}$ decays like $n^{-1}$. 
Although, the slope of the decay is dimension-independent, we note in Figure \ref{IS}  a dependence of the multiplicative constant on the dimension, i.e., the effective sample size is smaller in the case of $d\in\{3,4\}$ than in case of $d\in\{1,2\}$. 
This might be due the dependence of $\mathcal G$ on the dimension which effect the second and third derivatives at the modes -- which appear in our error bounds in the multiplicative constants. 

\begin{figure}[h]
    \centering
    \begin{minipage}{0.48\textwidth}
        \centering
        \includegraphics[width=0.9\linewidth]{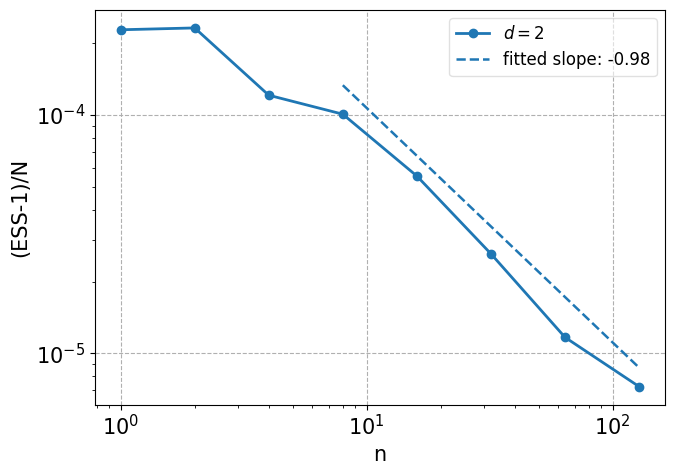}
        \caption{Convergence of $\left(\ESS-1\right)/N$ in prior-based importance sampling ($2D$ case)}
        \label{IS-prior}
    \end{minipage}
    \hfill
    \begin{minipage}{0.48\textwidth}
        \centering
        \includegraphics[width=0.9\linewidth]{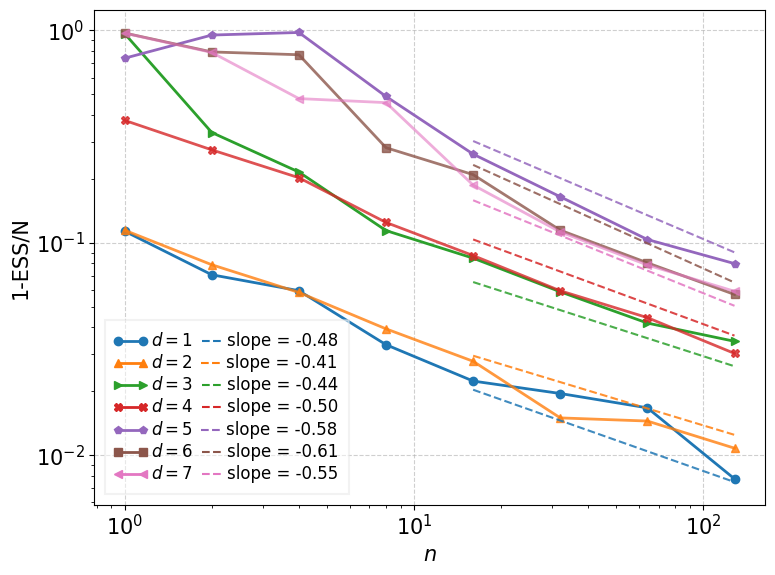}
        \caption{Convergence of $1-\ESS/N$ for Laplace-based independent Metropolis–Hastings algorithm}
        \label{IMH-ESS}
    \end{minipage}
\end{figure}

For comparison, we also plot $\left(\ESS-1\right)/N$ for prior-based importance sampling in the two-dimensional case with $N=10^5$ in Figure \ref{IS-prior}. The results show that $\left(\ESS-1\right)/N$ decays to 0 as $n$ increases, indicating that the ESS converges to $1$ as $n \to \infty$. This suggests that prior-based importance sampling becomes increasingly inefficient, whereas Laplace-based importance sampling becomes more efficient as $n\to\infty$.

\begin{figure}[H]
    \centering
    \begin{subfigure}[H]{0.4\textwidth}
        \includegraphics[width=\linewidth]{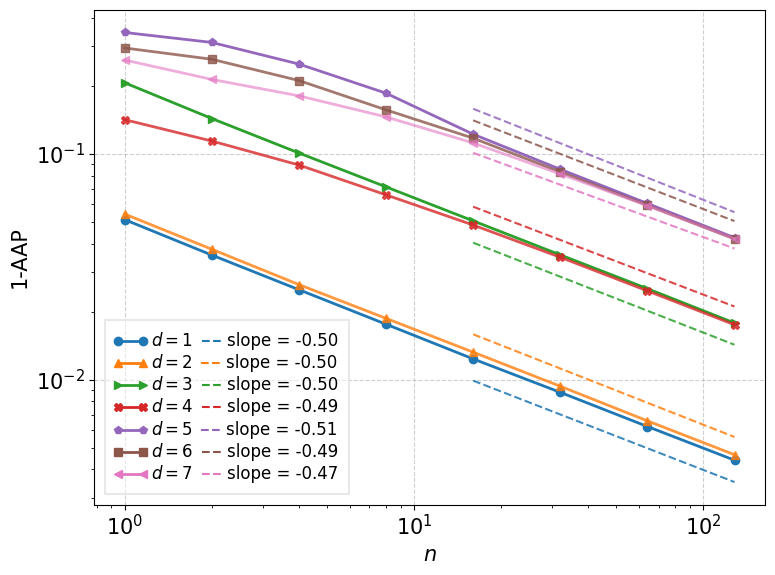}
        \caption{}
    \end{subfigure}
        \quad
    \begin{subfigure}[H]{0.4\textwidth}
        \includegraphics[width=\linewidth]{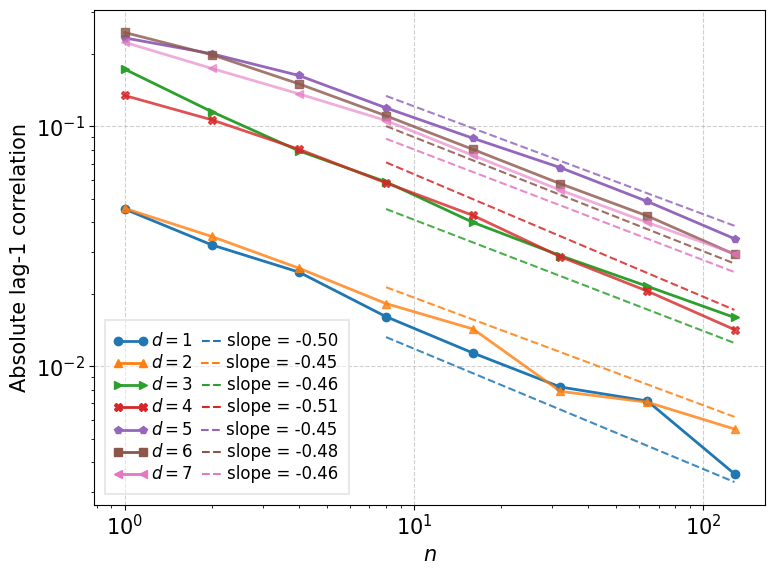}
        \caption{}
    \end{subfigure}
    \caption{(a): Convergence of the deviation of the average acceptance probability from 1; (b): Convergence of the absolute value of lag-one autocorrelation}
    \label{Lap-IMH}
\end{figure}

\paragraph{Laplace-based IMH:} We only consider the small noise limit here and examine how the average acceptance probability and the lag-one autocorrelation behave as $n$ increases. To this end, we run a single Markov chain with $10^6$ iterations to estimate the average acceptance probability. As shown in panel (a) of Figure \ref{Lap-IMH}, the average acceptance probability converges to one with a rate of $n^{-1/2}$ where the constant seems to depend on the dimension again. 
The lag one autocorrelation is also estimated from a single run of a Markov chain with $10^6$ iterations. We plot the absolute value of the lag-one autocorrelation in panel (b) of Figure \ref{Lap-IMH} and observe that it decays to zero at an approximate rate of $n^{-1/2}$. 
Also here the multiplicative constant seems to increase slightly with the dimension $d$ which is very likely caused by the dimension-dependence of $\mathcal G$.
Moreover, we also estimate the effective sample size (ESS) of Laplace-based IMH across different dimensions based on a single, long run of the Markov chain.
The ESS is given by number of iterations $N$ divided by the asymptotic variance $\widetilde \sigma_n^2(f)$. 
In Figure \ref{IMH-ESS} we observe that $1-\ESS/N$ decays also at a rate of $n^{-1/2}$. 
Due to the relation between spectral gap and ESS, i.e., $1 - \ESS/N \geq 1 - \frac{\mathrm{gap}_{\post}(M_n)}{2 \Vert f\Vert_{2,\post}}$, this numerical observation is in accordance with the statement of Theorem \ref{theo:sp} although its assumptions are not satisfied here.

\subsection{Bayesian neural network}
Neural networks are concatenations of simple functions consisting of affine layers and non-linear activation functions:
\[
    F(\xi) 
    =
    F_L \circ \cdots \circ F_1(\xi),
    \qquad
    \xi \in \R^N
\]
where each $F_k\colon \R^{N_{k-1}} \to \R^{N_k}$ is given by
\[
    F(\xi) 
    = 
    \phi\left(W_k  + b_k\right),
    \qquad
    W_k \in \R^{N_k \times N_{k-1}}, \quad b_k \in \R^{N_k},
\]
with chosen activation function $\phi: \mathbb{R} \to \mathbb{R}$ applied component-wise and $N_k$ denoting the number of nodes in the $k$th layer of the network. The values in the vectors $b_k$ are called \emph{biases} and the entries in the matrices $W_k$ are called \emph{weights}. In machine learning these biases and weights are the parameters to be trained given observational data of the unknown function to be approximated by the neural network.

In order to capture uncertainty in the model, i.e., uncertainty about the learned biases and weights, a classical approach is to consider \emph{Bayesian neural networks} where prior distributions for weights and biases are conditioned on observational data and the resulting posterior quantifies the limited knowledge about the neural network. 



\begin{figure}[H]
    \centering
    \begin{subfigure}[H]{0.45\textwidth}
        \includegraphics[width=\linewidth]{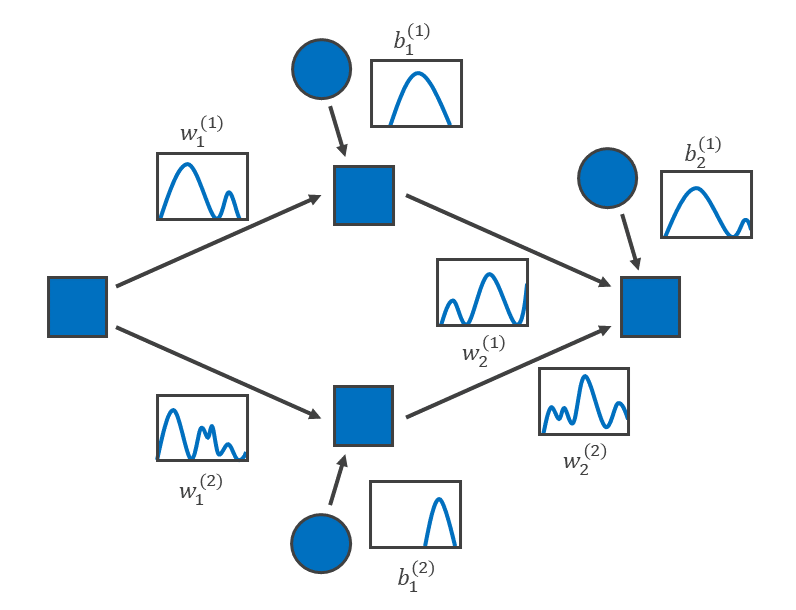}
        \caption{}
        \label{BNN}
    \end{subfigure}
        \quad
    \begin{subfigure}[H]{0.45\textwidth}
        \includegraphics[width=\linewidth]{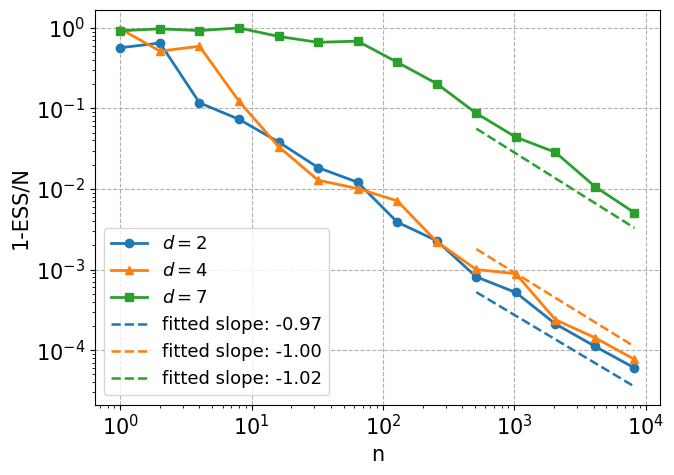}
        \caption{}
        \label{BNN_ESS}
    \end{subfigure}
    
    \caption{(a): Bayesian Neural Network with 
    $W_1 = \begin{bmatrix}
        w_1^{(1)} & w_1^{(2)}
    \end{bmatrix},
    b_1 = \left(b_1^{(1)},b_1^{(2)}\right),
    W_2 = \begin{bmatrix}
        w_2^{(1)}\\ w_2^{(2)}
    \end{bmatrix},$ and
    $b_1 = b_2^{(1)}$. (b): The convergence of effective sample size over the total sample size to 1 in Laplace-based importance sampling.}
\end{figure}

For simplicity, we consider a Bayesian neural network $F\colon \R \to \R$ with just a single hidden layer $L=2$ of two neurons, i.e., $N_0 = 1$, $N_1=2$, and $N_2=1$, cf. Figure \ref{BNN}, of the form
\[
    F(\xi) = F_2\left( \phi\left(F_1\left(\xi\right)\right) \right),
\]
where we choose as activation function the tanh function $\phi(t) = \frac{\mathrm e^{t} - \mathrm e^{-t}}{\mathrm e^{t} +\mathrm e^{-t}}$ and collect all parameters, i.e., all weights and biases, in a single vector
\[
    x = \left( w_1^{(1)},w_1^{(2)},b_1^{(1)},b_1^{(2)},w_2^{(1)},w_2^{(2)},b_2^{(1)} \right)^\top \in \mathbb{R}^7.
\]
The neural network, thus, becomes a parametrized mapping $F(\cdot; x)\colon \R \to \R$. 
We then assume a prior distribution $\pi_0 = \mathcal{N}\left( 0, 20 \cdot I_d\right)$ for the parameters $x \in \R^d$ and noisy observations of the neural network 
\[
    y_n 
    = 
    \mathcal G(X) + \eta_n,
    \qquad
    \mathcal G(x)
    =
    \begin{pmatrix}
        F(\xi_1;x) & \ldots & F(\xi_m;x)
    \end{pmatrix}^\top,
    \quad
    \eta_n \sim \mathcal{N}\left(0, n^{-1}\, \sigma^2\, I_m\right)
\]
at $m=100$ equidistant locations $\xi_1,\ldots,\xi_m \in [-5,5]$. Here, we choose $\sigma^2 = 10$ and assume as usual $X\sim \pi_0$ and $\eta_n$ to be independent. Thus, the resulting posterior $\post$ for the parameters $x$ is of the form \eqref{eq:definitionofmu} with
\[
    \logl_n := \frac{1}{20} \|y_n - \mathcal G(x)\|^2.
\]
In the experiment we generated the observed data by $y_n = \mathcal G(x_{\text{true}}) + \eta_n$ using a ground truth
\[
    x_{\text{true}}
    =
    \left( 1,-2,3,3,2,2,1 \right)^\top.
\]
We then consider three cases derived from fixing some of the parameters of the network to the true values and inferring only the $d\in\{2,4,7\}$ remaining parameters:
\begin{itemize}
\item 
2-dimensional case: we infer only the weights of the first layer, i.e., $w_1^{(1)}$ and $w_1^{(2)}$. 

\item 
4-dimensional case: we infer only all the weights, i.e., $w_1^{(1)}$, $w_1^{(2)}$, $w_2^{(1)}$, and $w_2^{(2)}$

\item 
7-dimensional case: we infer all seven parameters. 
\end{itemize}
Again, we apply Laplace-based importance sampling and Laplace-based IMH to sample from the multimodal posterior $\post$ and evaluate their performance for increasing $n$ in the various dimensions. In fact, for the cases $d\in\{2,4\}$ we have $K=2$ modes and for $d=7$ the posterior has $K=8$ modes.

\begin{figure}[h]
    \centering
    \begin{minipage}{0.48\textwidth}
        \centering
        \includegraphics[width=0.9\linewidth]{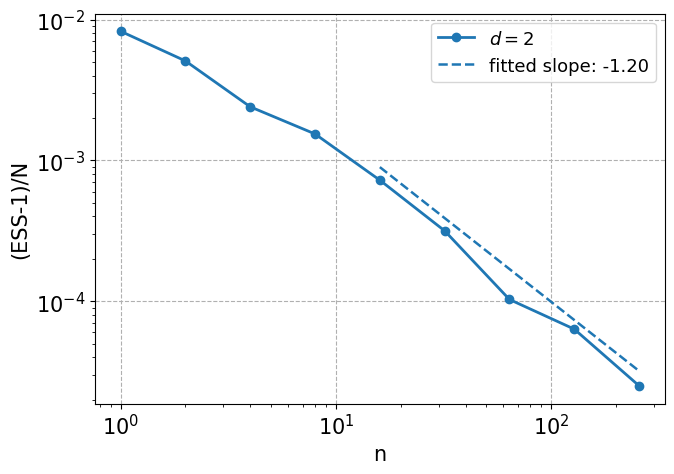}
        \caption{Convergence of $\left(1-\ESS\right)/N$ in prior-based importance sampling (2D case)}
        \label{BNN-IS-prior}
    \end{minipage}
    \hfill
    \begin{minipage}{0.48\textwidth}
        \centering
        \includegraphics[width=0.9\linewidth]{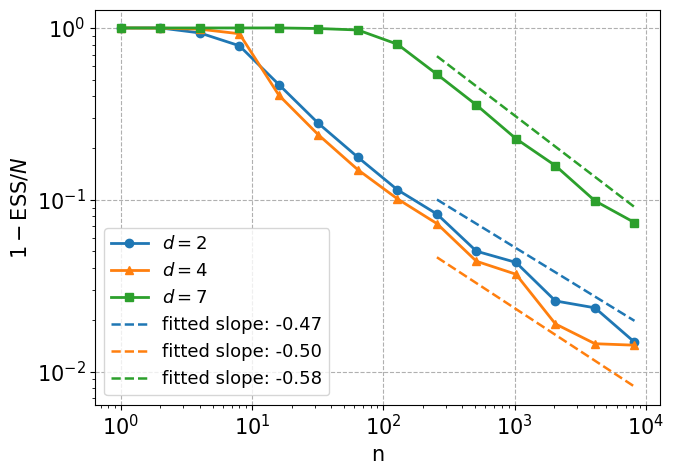}
        \caption{Convergence of $1-\ESS/N$ for Laplace-based independent Metropolis–Hastings algorithm}
        \label{BNN-IMH-ESS}
    \end{minipage}
\end{figure}

\paragraph{Results for Laplace-based importance sampling:} We use $N=10^3$ samples for the self-normalized Laplace-based importance sampling and study how the effective sample size evolves with increasing $n$. As shown in Figure \ref{BNN_ESS} the effective sample size converges to the total sample size $N$ with rate $n^{-1}$ in agreement with our theoretical findings. Again the multiplicative constant seems to depend on the dimension $d$, which is likely caused by the implicit dimension-dependence of $\mathcal G$ and its second and third derivatives at the modes. Note, that the lower-dimensional cases ($d=2$, $d=4$) are derived by fixing some of the seven parameters. Therefore, the norms of the second and third derivative of $\mathcal G$ at the modes is expected to grow with $d$ which can cause the observed decrease of ESS with $d$. 
Moreover, to compare with Laplace-based importance sampling, we plot $\left(\ESS-1\right)/N$ for prior-based importance sampling in the two-dimensional case with $N=10^5$. As $n$ increases, $\left(\ESS-1\right)/N$ decreases toward zero, indicating that the ESS converges to 1 as $n\to\infty$, as shown in Figure \ref{BNN-IS-prior}.

\begin{figure}[H]
    \centering
    \begin{subfigure}[H]{0.45\textwidth}
        \includegraphics[width=\linewidth]{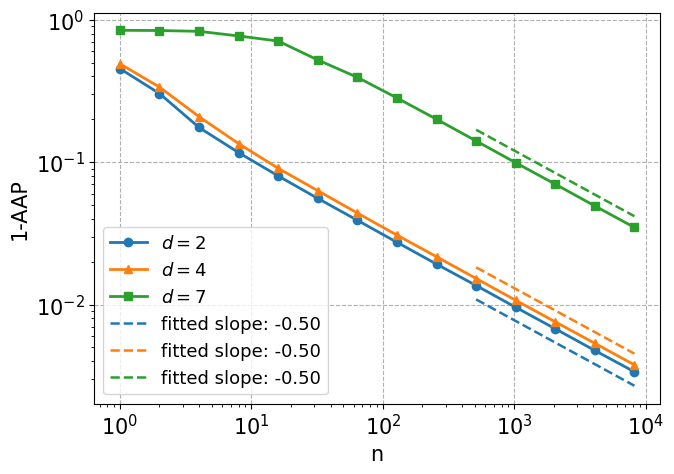}
        \caption{}
    \end{subfigure}
        \quad
    \begin{subfigure}[H]{0.45\textwidth}
        \includegraphics[width=\linewidth]{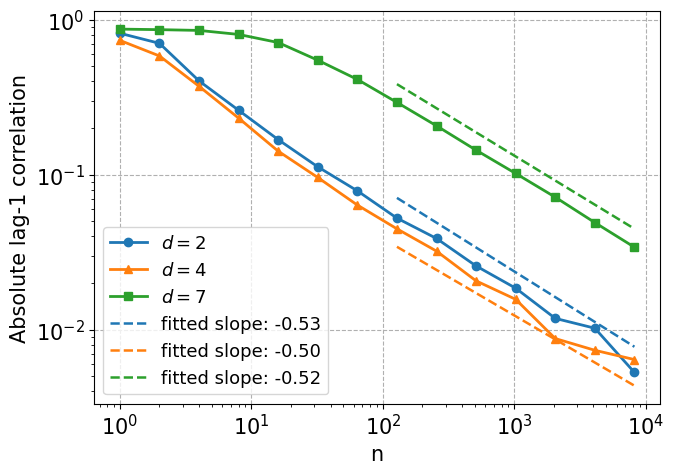}
        \caption{}
    \end{subfigure}
    
    \caption{Laplace-based independent Metropolis--Hastings algorithm: (a): Convergence of the deviation of the average acceptance probability from 1. (b): Convergence of the absolute value of lag-one autocorrelation.}
    
    \label{BNN_IMH}
\end{figure}

\paragraph{Results for Laplace-based-IMH:} As before we run a single Markov chain with $10^6$ iterations to estimate the average acceptance probability, the lag-one autocorrelation and the effective sample size. As shown in Figure \ref{BNN_IMH}, the average acceptance probability converges again to one and the lag-one autocorrelation again to zero each with rate $n^{-1/2}$ according to our theory.
Moreover, Figure \ref{BNN-IMH-ESS} shows that the effective sample size converges to the total number of iterations at rate $n^{-1/2}$.
Also for Laplace-based IMH we observe the same dimension dependence of the multiplicative constant as for Laplace-based SNIS.

\section{Conclusion} \label{sec:concl}
In this work we analysed a multi-modal Laplace approximation $\mathcal{L}^K_{\pi_n}$, i.e., a
weighted Gaussian mixture centred at the local minima of $I_n$, for sequences of concentrating
target measures $\pi_n$ with finitely many modes. We showed that the Hellinger distance between
$\pi_n$ and $\mathcal{L}^K_{\pi_n}$ decays as $n^{-1/2}$, which generalises known results for the
uni-modal case. Under the stronger tail condition of Assumption~3, we additionally obtained a
fully explicit, non-asymptotic bound. Its constants depend only polynomially on the dimension
$d$, and it also yields a bound on the $\chi^2$-divergence.

Building on these approximation results, we studied two Monte Carlo methods that use the Laplace
mixture as proposal. For self-normalised importance sampling, the effective sample size converges
in probability to the total sample size at a rate arbitrarily close to $n^{-1}$. Prior-based
importance sampling, in contrast, deteriorates at rate $n^{-d/2}$. For the independent
Metropolis--Hastings algorithm, the average acceptance probability converges to one and the
lag-one autocorrelation to zero, both at rate $n^{-1/2}$. Under stronger tail assumptions, the
spectral gap converges to one as well. Hence, both methods become more efficient as the target
concentrates rather than less, and the proposal does not need to be tuned. Our numerical
experiments confirm these rates. They
indicate that the same behaviour holds in the large data limit, which our theory does not cover.

Several questions remain open. 
For instance, our results on the effective sample size require all local minima of $I_n$ to
converge to global minimisers of a fixed potential $\Phi$, which excludes vanishing modes.
Similarly, the spectral gap result requires even $\delta > 1/2$ in Assumption~3, a condition that is
not always met in Bayesian inference. 
Without it, the independent Metropolis--Hastings algorithm may converge only subgeometrically, and studying this regime is a natural next step. Finally, our analysis assumes that all relevant modes are known. In practice, they have to be found by numerical optimisation, and the effect of missed modes or inexact Hessians on the sampling
efficiency deserves further study.


\appendix

\section{Proof of Hellinger convergence}\label{sec:proof}
Before proceeding to the proof of the main results, we state useful lemmas that allow us to keep as good as possible track on the dimensional dependence of our estimates


\subsection{Preliminary Lemmas}

This first lemma gives an explicit control on appearing moments of the standard multivariate normal distribution in $\R^d$.

\begin{lem}\label{lem:MomentLaplace} 
With $\Gamma$ denoting the Gamma function we have
for all $n \in \N$ and $p\in\{1,2\}$ that
    \begin{equation*}
        \int_{\R^d}\| x\|^{3p} \mathrm{e}^{-\frac{n}{2}\|x \|^2} \mathrm{d}x 
        = \frac{2^{\frac{3p}{2}}(2\pi)^{\frac{d}{2}}}{n^{\frac{3p+d}{2}}} \frac{\Gamma\big(\frac{3p+d}{2}\big)}{\Gamma\big(\frac{d}{2}\big)}
        \leq 2^{3p-2}
        \Big(\frac{2\pi}{n}\Big)^{\frac{d}{2}} \Big(\frac{d}{n}\Big)^{\frac{3p}{2}}
       .
    \end{equation*}
\end{lem}
\begin{proof}
The first equality is the $(3p)$-th moment of a multivariate normal distribution in $\R^d$, it follows from a standard change of variables to spherical coordinates and the representation of the Gamma function. For the upper bound we distinguish two cases: Firstly, for $p=2$ by the property $\Gamma(z+1)=z\cdot\Gamma(z)$, which holds for all $z>0$, we have
\begin{align*}
    \frac{\Gamma\big(\frac{d}{2}+3\big)}{\Gamma\big(\frac{d}{2}\big)}
    = \frac{d(d+2)(d+4)}{8} \leq \frac{d^3}{8}\Big(1+\frac{14}{d}\Big) \leq 2d^3,
\end{align*}
such that the inequality holds in that case. Secondly, for $p=1$, exploiting the aforementioned property and $\Gamma(z+1/2) \leq \sqrt{z} \cdot \Gamma(z)$, which holds for all $z>0$ (e.g. see \cite[Lemma~6]{novak}), we have
\begin{align*}
    \frac{\Gamma\big(\frac{d}{2}+\frac{3}{2}\big)}{\Gamma\big(\frac{d}{2}\big)}
    = \frac{\left( \frac{d}{2}+\frac{1}{2} \right)\Gamma\big(\frac{d}{2}+\frac{1}{2}\big)}{\Gamma\big(\frac{d}{2}\big)}
    \leq \frac{(d+1)\sqrt{d}}{2\sqrt{2}}  \leq \frac{d^{3/2}}{\sqrt{2}},
\end{align*}
such that the desired bound follows. 
\end{proof}

The next lemma gives an upper bound on Gaussian tails in arbitrary dimension.

\begin{lem}\label{lem:queue_dist_gauss} Let $\Sigma \in \R^{d\times d}$ be a symmetric definite positive matrix, $z\in \R^d$ and $r, c >0$. Then 
\begin{equation*}
     \int_{\left\|x - z \right\|^2_{\Sigma} > r }  \mathrm{e}^{-n c   \left\|x - z \right\|^2_{\Sigma}} \,\mathrm{d}x \le \mathrm{e}^{-\frac{n c r }{2}} 
     \mathrm{det}\big( \Sigma^{-1/2} \big)  
     \left(\frac{2\pi}{n c}\right)^{\frac{d}{2}}.
\end{equation*}
\end{lem}
\begin{proof}
A change of variable $x\mapsto \Sigma^{1/2}\left( x -z \right)$ 
with Jacobian of the inverse being $\Sigma^{-1/2} $ leads to
\begin{equation*}
     \begin{split}
         \int_{\left\|x - z \right\|^2_{\Sigma} > r } 
         & \mathrm{e}^{-n c   \left\|x - z \right\|^2_{\Sigma} }\, \mathrm{d}x = \int_{\|y\|^2 > r}  \mathrm{e}^{-n c \,\|y\|^2}\, \mathrm{det}\big(\Sigma^{-1/2} \big)  \mathrm{d}z
          = \mathrm{det}\big(\Sigma^{-1/2} \big) \left(\frac{\pi}{n c}\right)^{\frac{d}{2}} \P \left(\|Z\|^2 \ge r \right), 
    \end{split}
\end{equation*}
where $Z \sim \mathcal{N}(0,(2nc)^{-1} \mathrm{Id}_d)$.
By applying Markov's inequality, as $\P (\|Z\|^2 \ge r) \le \frac{\E [h(\|Z\|^2)]}{h(r)} $ with the increasing function $h(s) = \mathrm{e}^{\frac{n c s }{2}}$, and $\E \big[\mathrm{e}^{ \frac{n c}{2} \|Z\|^2}\big]=  2^{\frac{d}{2}}$ we obtain the desired estimate.
\end{proof}


Finally for the convenience of the reader we add a standard inequality without proof,  which is a few times used in the following in different settings. We have
\begin{equation}
\label{eq: Jensen_conseq}
\big \vert \sum_{j\in J} a_j \big\vert^p \leq \vert J \vert^{p-1} \sum_{j\in J} \vert a_j \vert^p 
\end{equation}
for any $p\geq1$, any finite index set $J$ with cardinality $\vert J \vert$ and arbitrary $(a_j)_{j\in J}\subset \R$.

\subsection{Main Lemma for proving \cref{theo:MainpropHelling2} and \cref{theo:MainpropHelling}}
\label{app:Mainlem}
\begin{lem}\label{lem:Mainlem} 
With the objects of the standing assumption, for any $n\in\N$ and $p\in\{1,2\}$, define
 \[
    \mathcal{K}_{n,p}
    :=
    \frac{2}{\widetilde Z_n^K} \int_{\R^d} \Big\vert \mathrm{e}^{-\frac{n}{p}\logpostd(x)}-\Big(\sum_{j=1}^K \pi_{n}^{(j)}\ \exp\big(-\frac{n}{2}\|x-x_n^{(j)}\|^2_{H_n^{(j)}}\big)\Big)^{\frac{1}{p}}  \Big\vert^p \mathrm{d}x,
 \]
 where $\widetilde Z_n^K$ is given as in \eqref{al: Z_tilde},
such that we have the following statements:
\begin{enumerate}[(a)]
\item 
\label{en: stat_a}
Let
Assumption~\ref{assum:Two} and Assumption~\ref{assum:coerc} with
$q=0$ and $\rho,M \in (0,\infty)$ as well as $R \leq \frac{\rho^6}{16M^2}$ and $\varepsilon\in (0,1)$ be satisfied. Then, for $p\in \{1,2\}$ and $n\in\N$, holds
 \begin{equation}\label{eq:Mainlem1}
  \mathcal{K}_{n,p}
  \leq 2^{\frac{d}{2} -2}  \bigg(\frac{2^{\frac{11}{2}} d^{\frac{3}{2}} M }{\rho^3 p }\bigg)^p n^{-\frac{p}{2}}  + \mathcal{O}\big(\mathrm{e}^{-\frac{n}{4} \min\{R,\varepsilon^2\}}\big).
\end{equation}
\item 
\label{en: stat_b}
Let \cref{assum:Two} and \cref{assum:Twobis}
with $q \in \left[\frac{1}{3}, \frac{1}{2}\right)$ and $\rho,M,\delta \in (0,\infty)$ as well as $R = \rho^2/ M^{\frac{2}{3}}$ be satisfied. 
Then, for $p\in \{1,2\}$ and $n\in\N$, holds
\begin{equation}\label{eq:Mainlem2}
   \mathcal{K}_{n,p}
  \leq
   \frac{\mathrm{e}}{4} \left(\frac{2^{4} d^{\frac{3}{2}} M }{\rho^3 p }\right)^p n^{-\frac{p}{2}} + 2^{p+\frac{d}{2}} \left( K (K-1)^{p-1}+1+(2\delta)^{-\frac{d}{2}} \right)\exp\Big(-\frac{\rho^2\min\{1,2\delta\}}{4M^{2/3}}\, n^{1-2q} \Big).
\end{equation}
\end{enumerate}
\end{lem}
\begin{proof}
Define 
\begin{equation} 
\label{eq: I_n^0_def}
    \mathcal{I}_n^{(0)} :=
    2^{p-1} \int_{\overline{V}_n^{R,q}} \exp(-n \logpostd (x)) \, \mathrm{d}x + 2^{p-1} \sum_{k=1}^K \pi_{n}^{(k)}  \int_{\overline{V}_n^{R,q}} \mathrm{exp}\big(-\frac{n}{2}\|x-x_n^{(k)}\|^2_{H_n^{(k)}} \big)\, \mathrm{d}x
\end{equation}
and note that by virtue of \eqref{eq: Jensen_conseq} we obtain
\begin{equation} 
 \label{eq: joint_I^0_n}
\int_{\overline{V}_n^{R,q}} \Big| ( \mathrm{exp}(-n\logpostd(x)))^{1/p}-\Big(\sum_{j=1}^K \pi_{n}^{(j)} \mathrm{exp}\Big(-\frac{n}{2}\|x-x_n^{(j)}\|^2_{H_n^{(j)}} \Big)\Big)^{1/p}\,  \Big|^p  \mathrm{d}x 
\leq \mathcal{I}_n^{(0)}.
\end{equation}
For $k=1,\dots,K$ let
\begin{align*}
    \mathcal{I}_{n}^{(k)}
    & :=
    \int_{\Vnq} \Big| ( \mathrm{exp}(-n\logpostd(x)))^{1/p}-\Big(\sum_{j=1}^K \pi_{n}^{(j)} \mathrm{exp}\Big(-\frac{n}{2}\|x-x_n^{(j)}\|^2_{H_n^{(j)}} \Big)\Big)^{1/p} \,\Big|^p \mathrm{d}x 
\end{align*}
such that, by decomposing $\R^d$ and \eqref{eq: joint_I^0_n}, we obtain
$
    \mathcal{K}_{n,p}
    \leq
    \frac{2}{\widetilde Z_n^K} \sum_{k=1}^K \mathcal{I}_n^{(k)}
    +
    \frac{2}{\widetilde Z_n^K} 
    \mathcal{I}_n^{(0)}.
$
We further estimate the terms of the appearing sum. For $k=1,\dots,K$, introducing
\begin{align*}
    \mathcal{I}^{(k),1}_n &:= \int_{\Vnq} \Big| ( \exp(-n \logpostd(x)))^{1/p} - \big( \pi_{n}^{(k)} \exp\big(-\frac{n}{2}\|x-x_n^{(k)}\|^2_{H_n^{(k)}}\big) \big)^{1/p} \, \Big|^p \mathrm{d}x\\
    \mathcal{I}^{(k),2}_n &:= \int_{\Vnq} \Big|( \pi_{n}^{(k)} \exp\big(-\frac{n}{2}\|x-x_n^{(k)}\|^2_{H_n^{(k)}}) \big)^{1/p}-\Big(\sum_{j=1}^K \pi_{n}^{(j)} \mathrm{exp}\Big(-\frac{n}{2}\left\|x-x_n^{(j)}\right\|^2_{H_n^{(j)}} \Big)\Big)^{1/p}  \,\Big|^p \mathrm{d}x
\end{align*}
applying \eqref{eq: Jensen_conseq} suitably yields
\begin{equation}\label{eq:Ikn}
\mathcal{I}^{(k)}_n \leq 2^{p-1} (\mathcal{I}^{(k),1}_n + \mathcal{I}^{(k),2}_n ).
\end{equation}
Furthermore
\begin{equation} \label{eq:Ikfacto}
    \begin{split}
        \mathcal{I}_{n}^{(k),1} & = \int_{\Vnq} \bigg| \bigg( \frac{\exp(-n \logpostd(x))}{\pi_n^{(k)} \exp\big(-\frac{n}{2}\|x-x_n^{(k)}\|^2_{H_n^{(k)}}\big)}  \bigg)^{1/p} -1 \bigg|^p \pi_{n}^{(k)} \exp\big({-\frac{n}{2}\|x-x_n^{(k)}\|^2_{H_n^{(k)}}}\big) \,\mathrm{d}x.
    \end{split}
\end{equation}
The second order Taylor expansion of $\logpostd$ gives for all $x \in \Vnq $ that
\begin{equation*}
        \logpostd(x) = \logpostd(x_n^{(k)}) + \langle \nabla \logpostd (x_n^{(k)} ), x-x_n^{(k)} \rangle  +\frac{1}{2} \|x-x_n^{(k)}\|^2_{H_n^{(k)}} + R_n^{(k)}(x) ,
\end{equation*}
with $R_n^{(k)} (x) := \frac{1}{6} D^3 \logpostd (\xi_x) (x - x_n^{(k)},x - x_n^{(k)},x - x_n^{(k)} ) $ for a vector $\xi_x \in \Vnq $. Since $x_n^{(k)}$ is a stationary point of $\logpostd$ we have $\nabla \logpostd(x_n^{(k)})=0$. Plugging the Taylor expansion into \eqref{eq:Ikfacto} yields
\begin{equation*}
        \mathcal{I}_{n}^{k,(1)}  = \pi_{n}^{(k)} \int_{\Vnq} \Big| \exp\Big(-\frac{n}{p} R^{(k)}_n(x)\Big) - 1 \Big|^p \exp\big(-\frac{n}{2}\|x-x_n^{(k)}\|^2_{H_n^{(k)}}\big)\, \mathrm{d}x,
\end{equation*}      
which, by using the inequality $|\exp(-t) - 1| 
\leq \exp(|t|)|t| $ that holds for all $t \in \R$, leads to
\begin{equation*}
       \mathcal{I}_{n}^{k,(1)}  \le \pi_{n}^{(k)} \int_{\Vnq} \exp(|n R^{(k)}_n(x)|) \,\Big|\frac{n}{p}R^{(k)}_n(x)\Big|^p \exp\big(-\frac{n}{2}\|x-x_n^{(k)}\|^2_{H_n^{(k)}}\big) \mathrm{d}x.
\end{equation*} 
Exploiting \cref{assum:Two}, in particular \eqref{en: M_bounded} and \eqref{en: unif_pos} (via Remark~\ref{rem: norm_ineq}), we obtain
    \begin{align}
    \notag
            \mathcal{I}_{n}^{k,(1)} & \leq \pi_{n}^{(k)} M^p \left(\frac{n}{p}\right)^{p} \int_{\Vnq}  \|x-x_n^{(k)} \|^{3p}  \exp\Big(-\frac{n}{2}\|x-x_n^{(k)}\|^2_{H_n^{(k)}}+n M \|x-x_n^{(k)} \|^3 \Big) \mathrm{d}x
            \\ 
    \label{al: I^1_start_for_further}        
            & \leq \pi_{n}^{(k)} M^p \left(\frac{n}{p}\right)^{p}  \int_{\Vnq}  \|x-x_n^{(k)} \|^{3p}  \exp\Big(-\frac{n}{2}\|x-x_n^{(k)}\|^2_{H_n^{(k)}}+ \frac{M n}{\rho^3} \|x-x_n^{(k)} \|^3_{H_n^{(k)}} \Big) \,\mathrm{d}x.
    \end{align}
Now we estimate $\mathcal{I}_n^{(k),2}$. For this observe that by \eqref{en: pairw_disj} of \cref{assum:Two} holds $\Vnq\subseteq \R^d \setminus V_n^{j,R,q}$ for any $j\in \{1,\dots,K\}$ with $j\neq k$.
Using this observation, applying \eqref{eq: Jensen_conseq} suitably and employing Lemma~\ref{lem:queue_dist_gauss} in the last estimate,
gives
\begin{align}
        \notag
        \mathcal{I}_n^{(k),2} 
        & = \int_{\Vnq} \Big|\sum_{j=1,\, j\neq k}^K (\pi_{n}^{(j)})^{1/p} \mathrm{exp}\Big(-\frac{n}{2p}\|x-x_n^{(j)}\|^2_{H_n^{(j)}} \Big)  \Big|^p \mathrm{d}x
        \\ 
        \notag
        & \leq (K-1)^{p-1} \sum_{j=1,\, j\neq k}^K \pi_{n}^{(j)} \int_{\Vnq} \mathrm{exp}\Big(-\frac{n}{2}\|x-x_n^{(j)}  \|^2_{H_n^{(j)}} \Big)  \mathrm{d}x
        \\ 
        \notag
        & \leq (K-1)^{p-1} \sum_{j=1,\,j\neq k}^K \pi_{n}^{(j)} \int_{\R^d \setminus V_n^{j,R,q}} \mathrm{exp}\Big(-\frac{n}{2}\|x-x_n^{(j)}\|^2_{H_n^{(j)}} \Big)  \mathrm{d}x \\
        \label{eq:Ikn22}
        & \leq 
        (K-1)^{p-1} \left(\frac{4\pi}{n }\right)^{d/2} \mathrm{e}^{-\frac{n^{(1-2q)} R }{4}} \sum_{j=1,\,j\neq k}^K \pi_{n}^{(j)} \,\mathrm{det}\big(( H_n^{(j)})^{-1/2} \big).
    \end{align}
In the following we further consider $\mathcal{I}_n^{(k),1}$, $\mathcal{I}_n^{(k),2}$ and $\mathcal{I}_n^{(0)}$ under the different assumptions imposed in statement \eqref{en: stat_a} and \eqref{en: stat_b}, that is, we treat both cases separately.

\paragraph{Addressing statement \eqref{en: stat_a}:} 
We further bound $\mathcal{I}_{n}^{(k),1}$ starting from \eqref{al: I^1_start_for_further}.
By the assumption that $q=0$ and $R \leq \frac{\rho^6}{16M^2}$ we get
\begin{align*}
       \mathcal{I}_{n}^{(k),1} 
       &  \leq \pi_{n}^{(k)} M^p  \Big(\frac{n}{p}\Big)^{p} \int_{\Vnz}  \|x-x_n^{(k)} \|^{3p}  \exp\Big(-\Big(\frac{n}{2}-\frac{n \sqrt{R} M}{\rho^3}  \Big)\|x-x_n^{(k)} \|_{H_n^{(k)}}^2\Big)\, \mathrm{d}x\\        
        & \leq \pi_{n}^{(k)} M^p \left(\frac{n}{p}\right)^{p} \int_{\Vnz}  \|x-x_n^{(k)} \|^{3p}  \exp\Big( -\frac{n}{4}\|x-x_n^{(k)}\|^2_{H_n^{(k)}} \Big)\, \mathrm{d}x
        \\ &  \leq \pi_{n}^{(k)} M^p \left(\frac{n}{p}\right)^{p}  \int_{\R^d}  \big\|\sqrt{2} (H_n^{(k)})^{-1/2}\varphi(x) \big\|^{3p}  \exp\big( -\frac{n}{2}\|\varphi(x)\|^2\big)\, \mathrm{d}x, 
\end{align*}
where $\varphi(x):= \frac{1}{\sqrt{2}}(H_n^{(k)})^{1/2}(x-x_n^{(k)} )$. Essentially, a change of variable yields 
\begin{align*}
        \mathcal{I}_{n}^{(k),1} & \leq \pi_{n}^{(k)} M^p \left(\frac{n}{p}\right)^{p}  2^{\frac{3p}{2}+ \frac{d}{2}}\vertiii{(H_n^{(k)})^{-1/2}}^{3p} \mathrm{det}\big( ( H_n^{(k)} )^{-1/2} \big)  \int_{\R^d}  \left\| y \right\|^{3p}  \mathrm{e}^{ -\frac{n}{2}\|y\|^2} \mathrm{d}y.
\end{align*}
Furthermore, Lemma~\ref{lem:MomentLaplace} and \eqref{en: unif_pos} of \cref{assum:Two} imply
\begin{align} 
        \mathcal{I}_{n}^{(k),1} & \leq \pi_{n}^{(k)} M^p \left(\frac{n}{p}\right)^{p} 2^{\frac{3p}{2}+ \frac{d}{2}} \,\vertiii{(H_n^k)^{-1/2}}^{3p}\, \mathrm{det}\big( ( H_n^{(k)} )^{-1/2} \big)  \, 
        2^{3p-2}\left(\frac{2 \pi}{n}\right)^{\frac{d}{2}}  \left(\frac{d}{n} \right)^{\frac{3p}{2}}  
        \notag
        \\ & 
        \label{eq:Ikn12}
        \leq \pi_{n}^{(k)} \left(\frac{2^{\frac{9}{2}} d^{\frac{3}{2}} M }{\rho^3 p }\right)^p  \left(\frac{2 \pi}{n}\right)^{\frac{d}{2}}\, 2^{\frac{d}{2}-2} \,\mathrm{det}\big( ( H_n^{(k)} )^{-1/2} \big)  n^{-\frac{p}{2}}.
\end{align}
Taking the definition of $\widetilde Z_n^K$ in \eqref{al: Z_tilde} into account and using the bounds of \eqref{eq:Ikn}, \eqref{eq:Ikn22} (with $q=0$) and \eqref{eq:Ikn12} yields 
\begin{align}
\notag
    \frac{2}{\widetilde Z_n^K} \sum_{k=1}^K \mathcal{I}_n^{(k)} & \leq \frac{2^p}{\widetilde Z_n^K} \sum_{k=1}^K (\mathcal{I}_n^{(k),1} + \mathcal{I}_n^{(k),2}) \\
\notag
 & \leq 2^{\frac{d}{2} -2}  \bigg(\frac{2^{\frac{11}{2}} d^{\frac{3}{2}} M }{\rho^3 p }\bigg)^p n^{-\frac{p}{2}} + 2^{p+\frac{d}{2}}(K-1)^{p-1}  \mathrm{e}^{-\frac{n R }{4}} \sum_{k=1}^K \frac{\sum_{j=1,\,j\neq k}^K \pi_{n}^{(j)} \mathrm{det}\big(( H_n^{(j)})^{-1/2} \big)}{\sum_{i=1}^K \pi_{n}^{(i)} \mathrm{det}\big(( H_n^{(i)})^{-1/2} \big)} \\ 
& 
\label{al: I_n^k-a-before-final}
\leq 
2^{\frac{d}{2} -2}  \bigg(\frac{2^{\frac{11}{2}} d^{\frac{3}{2}} M }{\rho^3 p }\bigg)^p n^{-\frac{p}{2}}+2^{p+\frac{d}{2}}K(K-1)^{p-1}  \mathrm{e}^{-\frac{n R }{4}}.
\end{align}

For deriving the desired bound of $\mathcal{K}_{n,p}$ with \eqref{al: I_n^k-a-before-final} we further estimate the ratio $2\mathcal{I}_{n}^{0}/\widetilde{Z}_n^K$. For this we consider the terms of $\mathcal{I}_{n}^{0}$ separately.
By the fact that $\overline{V}_n^{R,0} \subseteq\R^d\setminus V_{n}^{k,R,0}$ for any $k\in\{1,\dots,K\}$, Lemma~\ref{lem:queue_dist_gauss} and \eqref{al: Z_tilde} we obtain
\begin{align*}
    \sum_{k=1}^K \pi_{n}^{(k)}  \int_{\overline{V}_n^{R,0}} \mathrm{exp}\big(-\frac{n}{2}\|x-x_n^{(k)}\|^2_{H_n^{(k)}} \big)\,  \mathrm{d}x
  & \leq \sum_{k=1}^K \pi_{n}^{(k)}  \int_{\R^d\setminus V_{n}^{k,R,0}} \mathrm{exp}\big(-\frac{n}{2}\|x-x_n^{(k)}\|^2_{H_n^{(k)}} \big)\,  \mathrm{d}x \\
  & \leq 2^{\frac{d}{2}}\mathrm{e}^{-\frac{n R }{4}} \left(\frac{2\pi}{n}\right)^{\frac{d}{2}} \sum_{k=1}^K \pi_{n}^{(k)}
     \mathrm{det}\big( (H_n^{(k)})^{-1/2} \big) =  2^{\frac{d}{2}}\mathrm{e}^{-\frac{n R }{4}} \widetilde{Z}_n^{K},
\end{align*}
such that $\frac{2^{p}}{\widetilde{Z}_n^{K}}  \sum_{k=1}^K \pi_{n}^{(k)}  \int_{\overline{V}_n^{R,q}} \mathrm{exp}\big(-\frac{n}{2}\|x-x_n^{(k)}\|^2_{H_n^{(k)}} \big)\,  \mathrm{d}x \in \mathcal{O}(\mathrm{e}^{-\frac{n R }{4}})$. 

By Assumption~\ref{assum:coerc} for any $n\geq n^*$ and $x\in \overline{V}_n^{R,0}$ we have $\logpostd(x)\geq \varepsilon$, such that 
$g_n (x) := \mathbf{1}_{\overline{V}_n^{R,0}}(x) \exp(-n\logpostd(x)+n\varepsilon^2)$ satisfies 
$ g_n (x)  \le \exp(-n(1-\varepsilon)\varepsilon)$, which implies $\lim_{n\to\infty} g_n(x) = 0$. Furthermore, by \cref{assum:coerc} we obtain
\begin{equation*}
    g_n(x) \leq \exp(-n(1-\varepsilon)\logpostd(x)) \leq g(x)^{1-\varepsilon},
\end{equation*}
with $g^{1-\varepsilon}$ being Lebesgue integrable. Then, the dominated convergence theorem yields
\begin{equation*}
    \lim_{n\to\infty} \mathrm{e}^{n\varepsilon^2} \int_{ \overline{V}_n^{R,0}} \exp(-n \logpostd(x))\, \mathrm{d}x \leq \lim_{n\to\infty} \int_{\R^d} g_n(x) \mathrm{d}x 
    =  \int_{\R^d} \lim_{n\to\infty} g_n(x) \mathrm{d}x 
    = 0,
\end{equation*}
i.e., there is an $\widetilde{n}\in\N$ and $c_{R,\varepsilon}>0$, such that for all $n\geq \widetilde n $ holds $\int_{ \overline{V}_n^{R,0}} \exp(-n \logpostd(x))\, \mathrm{d}x\leq c_{R,\varepsilon} \mathrm{e}^{-n\varepsilon^2 }$.
By further exploiting Assumption~\ref{assum:coerc} we have $\widetilde Z_{n}^K \geq \big(\frac{2\pi}{n} \big)^{d/2} \varepsilon$, such that for any $n\geq \widetilde{n}$ holds
\[
 \frac{2^{p}}{\widetilde{Z}_n^K} \int_{\overline{V}_n^{R,0}} \exp(-n \logpostd (x)) \mathrm{d}x
 \leq \frac{2^p}{\varepsilon}\Big(\frac{n}{2\pi} \Big)^{d/2} c_{R,\varepsilon} \mathrm{e}^{-n\varepsilon^2 } \in \mathcal{O}(\mathrm{e}^{-n\frac{\varepsilon^2}{4} }),
\]
since $\sup_{n\in\N}\frac{n^{d/2}}{\mathrm{e}^{3n\varepsilon^2/4}} <\infty$.
Consequently, $\frac{2}{\widetilde Z_n^{K}}\mathcal{I}_n^{(0)} \in\mathcal{O}\big(\mathrm{e}^{-n\varepsilon^2/2 }+\mathrm{e}^{-\frac{n R }{4}}\big) \in \mathcal{O}(\mathrm{e}^{-\frac{n}{4} \min\{R,\varepsilon^2\}})$.

Now we conclude with the final estimate for statement \eqref{en: stat_a}. 
Combining \eqref{al: I_n^k-a-before-final} with $\frac{2}{\widetilde Z_n^K} \mathcal{I}_n^{(0)}\in \mathcal{O}(\mathrm{e}^{-\frac{n}{4} \min\{R,\varepsilon^2\}})$ gives 
\begin{equation} \label{eq:valueCrhoMeps}
    \mathcal{K}_{n,p}  \leq  2^{\frac{d}{2} -2}  \bigg(\frac{2^{\frac{11}{2}} d^{\frac{3}{2}} M }{\rho^3 p }\bigg)^p n^{-\frac{p}{2}}  + \mathcal{O}\big(\mathrm{e}^{-\frac{n}{4} \min\{R,\varepsilon^2\}}\big).
\end{equation}

\paragraph{Addressing statement \eqref{en: stat_b}:} 
We start again with bounding $\mathcal{I}_{n}^{(k),1}$ from \eqref{al: I^1_start_for_further}. For $q>1/3$ and 
$R = \rho^2/ M^{\frac{2}{3}}$ we get
\begin{align*}
         \mathcal{I}_{n}^{(k),1} 
        & \leq \pi_{n}^{(k)} M^p \left(\frac{n}{p}\right)^{p} \exp\Big(\frac{M R^{\frac{3}{2}}n^{1-3q}}{\rho^3} \Big) \int_{\Vnq}  \|x-x_n^{(k)} \|^{3p}  \exp\Big(-\frac{n}{2}\|x-x_n^{(k)}\|^2_{H_n^{(k)}} \Big) \,\mathrm{d}x\\        
        & \leq \mathrm{e}\, \pi_{n}^{(k)} M^p \left(\frac{n}{p}\right)^{p} \int_{\Vnq}  \|x-x_n^{(k)} \|^{3p}  \exp\Big(-\frac{n}{2}\|x-x_n^{(k)}\|^2_{H_n^{(k)}} \Big)\, \mathrm{d}x
        \\ &  \le \mathrm{e}\, \pi_{n}^{(k)} M^p \left(\frac{n}{p}\right)^{p}  \int_{\R^d}  \left\| (H_n^{(k)})^{-1/2} \varphi(x) \right\|^{3p}  \exp\Big(-\frac{n}{2}\|\varphi(x)\|^2\Big) \,\mathrm{d}x,
\end{align*}
where $\varphi(x) := (H_n^{(k)})^{1/2}(x-x_n^{(k)})$. 
Essentially, a change of variable yields 
\begin{equation*}
        \mathcal{I}_{n}^{(k),1}  \leq \mathrm{e}\, \pi_{n}^{(k)} M^p \left(\frac{n}{p}\right)^{p}  \vertiii{ (H_n^{(k)})^{-1/2}}^{3p} \mathrm{det}\big( ( H_n^{(k)} )^{-1/2} \big)  \int_{\R^d}  \left\| y \right\|^{3p}  \mathrm{e}^{ -\frac{n}{2}\|y\|^2} \mathrm{d}y.
\end{equation*}
Furthermore, Lemma~\ref{lem:MomentLaplace} and \eqref{en: unif_pos} of \cref{assum:Two} imply
\begin{align} 
        \notag
        \mathcal{I}_{n}^{(k),1} & \leq \mathrm{e}\, \pi_{n}^{(k)} M^p \left(\frac{n}{p}\right)^{p}   
        \vertiii{ (H_n^{(k)})^{-1/2}}^{3p} \mathrm{det}\big( ( H_n^{(k)} )^{-1/2} \big) \;2^{3p-2}\left(\frac{2 \pi}{n}\right)^{\frac{d}{2}} \left(\frac{d}{n} \right)^{\frac{3p}{2}} 
        \\ & 
        \label{eq:Ikn1}
        \leq \frac{\mathrm{e}}{4}\, \pi_{n}^{(k)} \left(\frac{2^{3} d^{\frac{3}{2}} M }{\rho^3 p }\right)^p  \left(\frac{2 \pi}{n}\right)^{\frac{d}{2}}\mathrm{det}\big( ( H_n^{(k)} )^{-1/2} \big)\, n^{-\frac{p}{2}}.
 \end{align}
Taking the definition of $\widetilde Z_n^K$ in \eqref{al: Z_tilde} into account and using the bounds of \eqref{eq:Ikn}, \eqref{eq:Ikn22} (with $R=\rho^2/M^{2/3}$) and \eqref{eq:Ikn1} yields 
\begin{align}
\notag
    & \frac{2}{\widetilde Z_n^K} \sum_{k=1}^K \mathcal{I}_n^{(k)}  \leq \frac{2^p}{\widetilde Z_n^K} \sum_{k=1}^K (\mathcal{I}_n^{(k),1} + \mathcal{I}_n^{(k),2}) \\
\notag
    & \leq \frac{2^p}{\widetilde Z_n^K} \sum_{k=1}^K \bigg[
    \frac{\mathrm{e}}{4}\, \pi_{n}^{(k)} \left(\frac{2^{3} d^{\frac{3}{2}} M }{\rho^3 p }\right)^p  \left(\frac{2 \pi}{n}\right)^{\frac{d}{2}}\mathrm{det}\big( ( H_n^{(k)} )^{-1/2} \big) n^{-\frac{p}{2}} \\
\notag
    & \qquad \qquad \qquad \qquad
    + (K-1)^{p-1} \left(\frac{4\pi}{n }\right)^{\frac{d}{2}} \exp\Big(-\frac{\rho^2 n^{(1-2q)}}{4M^{\frac{2}{3}}}\Big) \sum_{j=1,\,j\neq k}^K \pi_{n}^{(j)} \,\mathrm{det}\big(( H_n^{(j)})^{-1/2} \big)\bigg]\\
\notag   
& = \frac{\mathrm{e}}{4} \Big(\frac{2^{4} d^{\frac{3}{2}} M }{\rho^3 p }\Big)^p n^{-\frac{p}{2}} + 2^{p+\frac{d}{2}}(K-1)^{p-1} \exp\Big(-\frac{\rho^2 n^{(1-2q)}}{4M^{\frac{2}{3}}}\Big) \sum_{k=1}^K \frac{\sum_{j=1,j\neq k}^K \pi_{n}^{(j)} \mathrm{det}\big(( H_n^{(j)})^{-1/2} \big)}{\sum_{i=1}^K \pi_{n}^{(i)} \mathrm{det}\big(( H_n^{(i)})^{-1/2} \big)} \\
\label{al: I_n^k-b-before-final}
& \leq \frac{\mathrm{e}}{4} \Big(\frac{2^{4} d^{\frac{3}{2}} M }{\rho^3 p }\Big)^p n^{-\frac{p}{2}} + 2^{p+\frac{d}{2}} K (K-1)^{p-1} \exp\Big(-\frac{\rho^2 n^{(1-2q)}}{4M^{\frac{2}{3}}}\Big).
\end{align}   

Now we further estimate $\mathcal{I}_{n}^{(0)}$ with the aim to suitably control the ratio $\mathcal{I}_{n}^{(0)}/\widetilde Z_n^K$. By \cref{assum:Twobis}
the first term of $\mathcal{I}_{n}^{(0)}$ (cf. \eqref{eq: I_n^0_def}) satisfies
\begin{equation*}
            2^{p-1}\int_{\overline{V}_n^{R,q}} \exp(-n \logpostd (x))\, \mathrm{d}x 
            \leq 2^{p-1} \int_{\overline{V}_n^{R,q}} \sum_{k=1}^K  \pi_n^{(k)} \exp\big( -n \delta \|x - x_n^{(k)} \|_{H_n^{(k)}}^2 \big) \,\mathrm{d}x
\end{equation*}
such that
\begin{align*}
    \mathcal{I}_n^{(0)} & \le  2^{p-1} \int_{\overline{V}_n^{R,q}} \sum_{k=1}^K 
    \Big[\pi_{n}^{(k)} \mathrm{exp}\Big(-\frac{n}{2} \|x-x_n^{(k)}\|^2_{H_n^{(k)}} \Big) + \pi_n^{(k)} \exp\big(-n \delta   \|x - x_n^{(k)} \|^2_{H_n^{(k)}} \big) \Big] \, \mathrm{d}x 
    \\ & \le 2^{p-1} \sum_{k=1}^K \pi_{n}^{(k)}  \int_{\R^{d} \backslash \Vnq}  \mathrm{exp}\Big(-\frac{n}{2} \|x-x_n^{(k)}\|^2_{H_n^{(k)}} \Big) + \exp\big(-n \delta   \|x - x_n^{(k)} \|^2_{H_n^{(k)}} \big) \mathrm{d}x,
\end{align*} 
where the last inequality follows by the fact that $\overline{V}_n^{R,q} \subseteq \Vnq$ for any $k\in\{1,\dots,K\}$.
Applying Lemma~\ref{lem:queue_dist_gauss} twice leads to
\begin{align*}
    \mathcal{I}_n^{(0)} & 
    \leq 2^{p-1} \sum_{k=1}^K \pi_{n}^{(k)}  \mathrm{det}\big( ( H_n^{(k)} )^{-1/2} \big)  \bigg(\Big(\frac{4\pi}{n}\Big)^{\frac{d}{2}} \mathrm{e}^{-\frac{R}{4}\,n^{1-2q}} +  \Big(\frac{2\pi}{n \delta}\Big)^{\frac{d}{2}} \mathrm{e}^{-\frac{\delta R}{2}\,n^{1-2q}} \bigg)
    \\ & \le 2^{p-1} \left(  4^{\frac{d}{2}} \mathrm{e}^{-\frac{R}{4}\,n^{1-2q}} + 2^{\frac{d}{2}} \delta^{-\frac{d}{2}} \mathrm{e}^{-\frac{\delta R}{2}\,n^{1-2q}} \right) \left(\frac{\pi}{n}\right)^{\frac{d}{2}}  \sum_{k=1}^K \pi_{n}^{(k)}  \mathrm{det}\big( ( H_n^{(k)} )^{-1/2} \big) 
    \\ & \le 2^{p-1+\frac{d}{2}}  \mathrm{e}^{-\frac{R}{4}\min\{1,2\delta\}\, n^{1-2q}  } 
    \left(  2^{\frac{d}{2}}  + \delta^{-\frac{d}{2}} \right)
    \left(\frac{\pi}{n}\right)^{\frac{d}{2}}  \sum_{k=1}^K \pi_{n}^{(k)}  \mathrm{det}\big( ( H_n^{(k)} )^{-1/2} \big).
\end{align*}
Plugging \eqref{al: Z_tilde} and $R = \rho^2/M^{2/3}$ in, gives
\begin{equation*}
\frac{2}{\widetilde Z_n^K} 
    \mathcal{I}_n^{(0)}
\leq 2^{p} \left(  2^{\frac{d}{2}}  + \delta^{-\frac{d}{2}} \right)
\exp\Big(-\frac{\rho^2\min\{1,2\delta\}}{4M^{2/3}}\, n^{1-2q} \Big),
\end{equation*}
such that we finally obtain with \eqref{al: I_n^k-b-before-final} that
\begin{align*}
    \mathcal{K}_{n,p}
    & \le 
   \frac{\mathrm{e}}{4} \left(\frac{2^{4} d^{\frac{3}{2}} M }{\rho^3 p }\right)^p n^{-\frac{p}{2}} + 2^{p+\frac{d}{2}} ( K (K-1)^{p-1}+1+(2\delta)^{-\frac{d}{2}} )\exp\Big(-\frac{\rho^2\min\{1,2\delta\}}{4M^{2/3}}\, n^{1-2q} \Big),
\end{align*}
which finishes the proof.
\end{proof}

\subsection{Proof of \cref{theo:MainpropHelling2} and \cref{theo:MainpropHelling}}
\label{sec: proof_main_Helling}
Define
\begin{align*}
    \widetilde{\mathcal{K}}_{n,1}
    & := 2 Z_n \big( (\widetilde Z_n^K)^{-1/2} - Z_n^{-1/2} \big)^2,\\
    \widetilde{\mathcal{K}}_{n,2} & := \frac{2}{\widetilde Z_n^K} \int_{\R^d} \Big \vert \mathrm{e}^{-\frac{n}{2}\logpostd(x)}-\Big(\sum_{j=1}^K\pi_{n}^{(j)} \,\mathrm{exp}\big(-\frac{n}{2}\|x-x_n^{(j)}\|^2_{H_n^{(j)}} \big)\Big)^{\frac{1}{2}}  \Big\vert^2 \mathrm{d}x
    \end{align*}
and note that by suitable application of \eqref{eq: Jensen_conseq} we have
\begin{equation*}
        d_{\mathrm{H}}^2( \lap ,\post ) =  \int_{\R^d} \left[\sqrt{\post(x)}-\sqrt{\lap(x)}  \right]^2 \mathrm{d}x
        \le  \widetilde{\mathcal{K}}_{n,1} +  \widetilde{\mathcal{K}}_{n,2}.
\end{equation*}
Moreover, by the definition of $\mathcal{K}_{n,1}, \mathcal{K}_{n,2}$ from Lemma~\ref{lem:Mainlem} we observe $\widetilde{\mathcal{K}}_{n,2} = \mathcal{K}_{n,2}$ and
\begin{align*}
        \widetilde{\mathcal{K}}_{n,1} &  
        = \frac{2 Z_n ( (\widetilde Z^K_n)^{-1}  -  Z_n^{-1} )^2 }{ \big( (\widetilde Z_n^K)^{-1/2} + Z_n^{-1/2} \big)^2}
        \leq \frac{2 Z_n ( (\widetilde Z^K_n)^{-1}  -  Z_n^{-1} )^2 }{ \big( Z_n^{-1/2} \big)^2}
        = \frac{2  ( Z_n - \widetilde Z_n^K  )^2 }{ (\widetilde Z_n^K) ^2}
        \leq \frac{\mathcal{K}_{n,1}^2}{2}.
\end{align*}
\begin{proof}[Proof of \cref{theo:MainpropHelling2}]

Applying statement \eqref{en: stat_a} of Lemma~\ref{lem:Mainlem} with $p=1$ and $p=2$, if necessary in combination with \eqref{eq: Jensen_conseq}, leads to
\begin{equation*}
    d_{\mathrm{H}}^2( \lap ,\post )
    \leq \frac{2^{11}d^3 M^2}{\rho^6} (2^{d-4}+2^{\frac{d}{2}-4}) n^{-1} + \mathcal{O}\big(\mathrm{e}^{-\frac{n}{4} \min\{R,\varepsilon^2\}}\big),
\end{equation*}
such that $ d_{\mathrm{H}}( \lap ,\post ) \in \mathcal{O}(n^{-\frac{1}{2}})$.
\end{proof}

\begin{proof}[Proof of \cref{theo:MainpropHelling}]
Applying statement \eqref{en: stat_b} of Lemma~\ref{lem:Mainlem} with $p=1$ and $p=2$ leads to
\begin{align*}
& d_{\mathrm{H}}^2( \lap ,\post ) 
 \leq \frac{\mathrm{e}2^4 d^3 M^2}{n \rho^6} (\mathrm{e} + 1) + 4\exp\Big(-\frac{\rho^2 \min\{1,2\delta\}}{2M^{\frac{2}{3}}} n^{(1-2q)}\Big)
    [2^d( K+1+(2\delta)^{-\frac{d}{2}})^2]\\
 & \qquad \qquad \qquad \quad +4\exp\Big(-\frac{\rho^2 \min\{1,2\delta\}}{4M^{\frac{2}{3}}} n^{(1-2q)}\Big) [2^{\frac{d}{2}}(K^2-K+1+(2\delta)^{-\frac{d}{2}})]\\
 & \leq \frac{\mathrm{e}2^4 d^3 M^2}{n \rho^6} (\mathrm{e} + 1) 
 + 2^{d+2}\exp\Big(-\frac{\rho^2 \min\{1,2\delta\}}{4M^{\frac{2}{3}}} n^{(1-2q)}\Big)
 [2(K+1)^2 + 2(2\delta)^{-d} + K^2 + (2\delta)^{-\frac{d}{2}}]\\
  & \leq \frac{\mathrm{e}2^4 d^3 M^2}{n \rho^6} (\mathrm{e} + 1) 
 + 12\cdot2^{d}\exp\Big(-\frac{\rho^2 \min\{1,2\delta\}}{4M^{\frac{2}{3}}} n^{(1-2q)}\Big)
 [(K+1)^2 + \max\{1,(2\delta)^{-d}\}],
 \end{align*}
where also \eqref{eq: Jensen_conseq} and $2(2\delta)^{-d} + (2\delta)^{-\frac{d}{2}}\leq 3\max\{1,(2\delta)^{-d}\}$ has been applied. Consequently, taking the square root and simple estimates give the desired bound
\begin{equation*}
d_{\mathrm{H}}( \lap ,\post ) 
\leq \frac{13 d^{\frac{3}{2}}M}{n^{\frac{1}{2}} \rho^3} 
+ 2^{\frac{d}{2}+2}  \exp\Big(-\frac{\rho^2 \min\{1,2\delta\}}{8M^{\frac{2}{3}}} n^{(1-2q)}\Big) [(K+1) + \max\{1,(2\delta)^{-\frac{d}{2}}\}].
\qedhere
\end{equation*}
\end{proof}
\subsection{Proof of \Cref{cor: accesable_statement}}
\label{subsec: proof_coro}
We start with a purely analytical auxiliary estimate.
\begin{lem} \label{lem: aux_for_acc_stat}
For $C>0$, $q\in[0,1/2)$ with $D\geq \frac{1+\mathrm{e}^{-1}}{1-2q}$ and for 
any $n\in\N$ with $n\geq \Big( \frac{D}{C} \log\big(\frac{D}{C}\big) \Big)^{1-2q}$
we have $\exp(-Cn^{1-2q}) \leq n^{-1}$.
\end{lem}
\begin{proof}
    Define $f\colon (0,\infty) \to (0,\infty)$ by $f(x):=\frac{\log x}{x^{1-2q}}$ and note that
    $\exp(-Cn^{1-2q}) \leq n^{-1} $ iff $f(n)\leq C$. It is easy to see that $f$ is maximized at $n^*:=\mathrm{e}^{\frac{1}{1-2q}}$ with $f(n^*)=((1-2q)\mathrm{e})^{-1}$. 
    We distinguish two cases: Firstly, if $C>f(n^*)$, then trivially $f(n)\leq f(n^*) \leq C$ for any $n\in\N$, such that the claim is true. Secondly, assume that $0<C\leq f(n^*)$.
    Here it is useful to observe that $f$ is strictly decreasing for $x>n^*$ and
    \[
    n\geq \Big( \frac{D}{C} \log\big(\frac{D}{C}\big) \Big)^{1-2q}
     \geq \Big( (\mathrm{e}+1)\log(\mathrm{e}+1) \Big)^{1-2q} > n^*,
    \]
    where we used that $\frac{D}{C}\geq \mathrm{e}+1$ as well as the fact that $x\mapsto x\log x$ is strictly increasing for $x>\mathrm{e}^{-1}$. Consequently, by $n>n^*$ and the monotonicity we have
    \begin{align*}
        f(n) & \leq 
        \frac{\log\Big[\big(\frac{D}{C}\big)^{1-2q} \big[\log \big(\frac{D}{C}\big)\big]^{1-2q}\Big]}{\frac{D}{C}\log\big(\frac{D}{C} \big)}
        = \frac{C}{D} \Big( \frac{1}{1-2q}+f\Big( \big[ \log \big( \frac{D}{C}\big) \big]^{1-2q}\Big) \Big)\\
        & \leq \frac{C}{D} \Big( \frac{1}{1-2q}+f(n^*) \Big) = \frac{C(1+\mathrm{e}^{-1})}{D(1-2q)} \leq C,
    \end{align*}
    where the last inequality follows by $D\geq  \frac{1+\mathrm{e}^{-1}}{1-2q}$.
\end{proof}
\begin{proof}[Proof of \cref{cor: accesable_statement}]
Note that \eqref{eq: conseq_of_thm} follows immediately from Theorem~\ref{theo:MainpropHelling2} and the fact that $\delta\geq1/2$.
We turn to the proof of the simplified estimate when $n$ is sufficiently large.
By Lemma~\ref{lem: aux_for_acc_stat} with $C=\frac{\rho^2}{8M^{2/3}}$, $D=\frac{2}{1-2q}$, and taking the square root, we obtain with $n_1:= \Big [\frac{16 M^{2/3}}{\rho^2(1-2q)} \log \big( \frac{16 M^{2/3}}{\rho^2(1-2q)} \big)\Big]^{1-2q}$, that
\begin{equation}
\label{eq: exp_by_poly}
\exp\Big(-\frac{\rho^2 n^{1-2q}}{16 M^{2/3}}\Big) \leq n^{-\frac{1}{2}} \qquad \text{for any $n\geq n_1$.}
\end{equation}
Furthermore, for $n_2 
:= \Big[\frac{16 M^{2/3}}{\rho^2} \log\big( \frac{2^{\frac{d}{2}+2}(K+2)\rho^3}{M} \big) \Big]^{\frac{1}{1-2q}}
\geq \Big[\frac{16 M^{2/3}}{\rho^2} \log\big( \frac{2^{\frac{d}{2}+2}(K+2)\rho^3}{d^{3/2}M} \big) \Big]^{\frac{1}{1-2q}}
$ we have
\begin{equation}
\label{eq: exp_by_const}
    \exp\Big(-\frac{\rho^2 n^{1-2q}}{16 M^{2/3}}\Big) \leq \frac{d^{\frac{3}{2}}M}{2^{\frac{d}{2}+2}(K+2)\rho^3}
    \qquad 
    \text{for any $n\geq n_2$.}
\end{equation}
Consequently, for any $n\in\N$ with $n\geq \max\{n_1,n_2\}$ we obtain
\begin{align*}
2^{\frac{d}{2}+2} (K+2)  \exp\Big(-\frac{\rho^2 n^{(1-2q)} }{8M^{\frac{2}{3}}}\Big)
\underset{\eqref{eq: exp_by_poly}}{\leq} 
\frac{2^{\frac{d}{2}+2} (K+2) }{n^{\frac{1}{2}}} 
\exp\Big(-\frac{\rho^2 n^{(1-2q)} }{16M^{\frac{2}{3}}}\Big)
\underset{\eqref{eq: exp_by_const}}{\leq}
\frac{M d^{\frac{3}{2}}}{n^{\frac{1}{2}}\rho^3},
\end{align*}
which in combination with \eqref{eq: conseq_of_thm} proves the desired statement.
\end{proof}
\subsection{Proof of \Cref{theo:PropChiSq}}
\label{sec:proofChiSq}
%
%
%
It is sufficient to verify that $\frac{\post(x)}{\lap(x)} \le 2\mathrm{e}$ for any $x\in\R^d$ and any $n\in\N$ that satisfies \eqref{eq: restr_chi_squ}. Then, the rest follows by \eqref{eq:chi2_hellinger} and \eqref{eq: conseq_of_thm}.

Firstly, we prove that $\frac{\post(x)}{\lap(x)} \le  \frac{\mathrm{e} \widetilde Z_n^K}{Z_n}$ for any $x\in\R^d$ and any $n\in\N$. For this we rely on the `ratio function' $\unimpw$, introduced in \eqref{weight_LM}, that is chosen in such a way that $\frac{\post(x)}{\lap(x)}=\frac{\widetilde Z_n^K}{Z_n} \unimpw(x)$. 
For $x\in \bigcup_{j=1}^K V_n^{j,R,q}$ there is an index $k\in\{1,\dots,K\}$ such that $x\in V_n^{k,R,q}$ and 
    \[
        \unimpw(x) = \frac{\mathrm{e}^{-n\logpostd(x)}}{\pi_n^{(k)} \exp{(- \frac{n}{2}\|x-x_n^{(k)}\|^2_{H_n^{(k)}})}\Big(1+ \sum_{j=1,\,j\neq k}^K \frac{\pi_n^{(j)}}{\pi_n^{(k)}} \exp(-\frac{n}{2}\|x-x_n^{(j)}\|^2_{H_n^{(j)}}+\frac{n}{2}\|x-x_n^{(k)}\|^2_{H_n^{(k)}}) \Big)},
    \]
    which implies
    \begin{equation*}
        \unimpw(x) \le \exp\left( -n \left(\logpostd(x)-I_n(x_n^{(k)})-\frac{1}{2}\|x-x_n^{(k)}\|^2_{H_n^{(k)}}\right)\right),
    \end{equation*}
    where we also used the definition of $\pi_n^{(k)}$. By the fact that $x\in\Vnq$
    and a second order Taylor expansion of $\logpostd$ we have
    \begin{equation*}
    \begin{split}
        \logpostd(x) & = \logpostd(x_n^{(k)}) + \big\langle \nabla \logpostd (x_n^{(k)} ), x-x_n^{(k)} \big\rangle  +\frac{1}{2} \|x-x_n^{(k)}\|^2_{H_n^{(k)}} + R_n^{(k)}(x) ,
    \end{split}
\end{equation*}
with $R_n^{(k)} (x) = \frac{1}{6} D^3 \logpostd (\xi_x) (x - x_n^{(k)},x - x_n^{[k]},x - x_n^{(k)} ) $ for some $\xi_x \in \Vnq $. 
Exploiting $\nabla \logpostd(x_n^{(k)})=0$, which holds by the fact that 
$x_n^{(k)}$ is a stationary point, as well as the Taylor expansion in combination with 
\cref{assum:Two} \eqref{en: M_bounded} and Remark~\ref{rem: norm_ineq}, we obtain
    \begin{equation*}
        \unimpw(x) \le  e^{n|R^{(k)}_n(x)|} \leq  \exp\Big( \frac{Mn}{\rho^3}  \|x- x_n^{(k)}\|_{H_n^{(k)}}^3 \Big).
    \end{equation*}
    By the definition of $\Vnq$, $R = \rho^2/ M^{\frac{2}{3}}$ and $x\in \Vnq$ we obtain
        \begin{equation}
        \label{eq:important_bnd_spectral_gap}
                    \unimpw(x)  \le  \mathrm{e}^{n^{1-3q}} \leq \mathrm{e},
        \end{equation}
    where the last inequality follows by $q \in [1/3,1/2)$.
For  $x\in\overline{V}_n^{R,q}$, by \cref{assum:Twobis},  we can conclude     
\begin{align*}
\unimpw (x) &\le     \frac{\sum_{k=1}^K \pi_n^{(k)} \exp\big(-n \delta  \|x - x_n^{(k)} \|^2_{H_n^{(k)}}\big)}{\sum_{j=1}^K \pi_n^{(j)}  \exp\big(-\frac{n}{2} \|x - x_n^{(j)} \|^2_{H_n^{(j)}} \big)}
\leq 1,
\end{align*}
where the last inequality
follows by $\delta\geq 1/2$ which implies
$\exp\big(-n \delta  \|x - x_n^{(j)} \|^2_{H_n^{(j)}}\big)
\leq \exp\big(-\frac{n}{2}  \|x - x_n^{(j)} \|^2_{H_n^{(j)}}\big)$ for any $j\in\{1,\dots,K\}$.
Consequently we proved that for any $x\in\R^d$ holds $\frac{\post(x)}{\lap(x)} \le  \frac{\mathrm{e} \widetilde Z_n^K}{Z_n}$.
With $\mathcal{K}_{n,1}$ defined in Lemma~\ref{lem:Mainlem} holds
\begin{align*}
   & \Big\vert \frac{Z_n}{\widetilde{Z}_n^K}-1 \Big\vert  
    = \frac{1}{\widetilde{Z}_n^K} \vert \widetilde{Z}_n^K-Z_n\vert
   \leq \frac{1}{\widetilde{Z}_n^K} \int_{\R^d} \Big\vert \mathrm{e}^{-n\logpostd(x)}-\sum_{j=1}^K \pi_{n}^{(j)}\ \exp\big(-\frac{n}{2}\|x-x_n^{(j)}\|^2_{H_n^{(j)}}\big)  \Big\vert \mathrm{d}x 
   = \frac{\mathcal{K}_{n,1}}{2},
\end{align*}
such that
\begin{equation}
\label{al:cite_for_spectral_gap}
1- \frac{\mathcal{K}_{n,1}}{2}\leq \frac{Z_n}{\widetilde{Z}_n^K} 
\qquad \text{implies} \qquad
\frac{\widetilde{Z}_n^K}{Z_n}\leq (1- \frac{\mathcal{K}_{n,1}}{2})^{-1} = \frac{2}{2-\mathcal{K}_{n,1}}.
\end{equation}
For $n$
satisfying \eqref{eq: restr_chi_squ} it follows
\begin{align*}
& \frac{4\mathrm{e}\, d^{3/2}M}{\rho^3 n^{1/2}} \leq \frac{1}{2} \qquad \text{and}\qquad
 2^{1+\frac{d}{2}}(K+2) \exp\Big(-\frac{\rho^2}{4M^{2/3}} n^{1-2q} \Big) \leq \frac{1}{2}
\end{align*}
such that by \eqref{en: stat_b} of Lemma~\ref{lem:Mainlem} (with $\delta\geq1/2$) 
we have $\mathcal{K}_{n,1}\leq 1$. Consequently (for those $n$) holds 
$\frac{\widetilde{Z}_n^K}{Z_n}\leq 2$ and $\frac{\post(x)}{\lap(x)} \le  2\mathrm{e}$ as desired.

\section{Proofs of Laplace-based importance sampling}
We provide proofs that have been postponed in Section~\ref{sec:LAIS}.

\subsection{Proof of Proposition \ref{prop: GnN}.}\label{sec:GnN_proof}
\begin{proof}
For any $k\in\{1,\dots,K\}$ and $i\in \{1,\dots,N\}$
observe that
    \[
        A_{i,n}
        :=
        \Big\{ \min_{j=1,\ldots,K}
         \| X_{i,n} - x_n^{(k)} \|^2_{H_n^{(k)}} > Rn^{-2q} \Big\} =\Big\{ X_{i,n} \in\bigcap_{j=1}^K 
         (V_{n}^{j,R,q})^c \Big\} \subseteq \{X_{i,n}\in (V_{n}^{k,R,q})^c \}.
   \]
     Consequently, we obtain
    \begin{align*}
        \mathbb{P}( A_{i,n})
        & 
         = \int_{\bigcap_{k=1}^K 
         (\Vnq)^c} \frac{1}{\widetilde{Z}_n^K} \sum_{k=1}^K \pi_n^{(k)} \mathrm{e}^{-\frac{n}{2} \|x - x_n^{(k)}\|^2_{H_n^{(k)}} } {\rm d}x 
        = \sum_{k=1}^K \frac{\pi_n^{(k)}}{\widetilde{Z}_n^K} \int_{{\bigcap_{j=1}^K 
         (V_{n}^{j,R,q})^c}} 
          \mathrm{e}^{-\frac{n}{2} \|x - x_n^{(k)}\|^2_{H_n^{(k)}} } {\rm d}x 
    \end{align*}
      and, by \cref{lem:queue_dist_gauss} we have
    \begin{align*}
        \int_{{\bigcap_{j=1}^K 
         (V_{n}^{j,R,q})^c}} \mathrm{e}^{ -\frac{n}{2} \|x - x_n^{(k)}\|^2_{H_n^{(k)}} }\, \mathrm d x 
         \leq \int_{(\Vnq)^c} \mathrm{e}^{-\frac{n}{2} \|x - x_n^{(k)}\|^2_{H_n^{(k)}} } {\rm d}x 
        \leq \mathrm e^{-\frac{Rn^{(1-2q)}}{4}} \Big( \frac{4\pi}{n} \Big)^{\frac{d}{2}} \det\big((H_n^{(k)})^{-1}\big)^{1/2}.
    \end{align*}
    By $\widetilde Z_n^K  = \big(\frac{2 \pi}{n} \big)^{d/2}\sum_{j=1}^K  \det\big( \big(H_n^{(j)}\big)^{-1/2} \big) \pi_{n}^{(j)}$ this gives
    \begin{align*}
        \mathbb{P}(A_{i,n}) 
        &\leq \sum_{k=1}^K \frac{\pi_n^{(k)}}{\widetilde{Z}_n^K} \mathrm{e}^{-\frac{Rn^{(1-2q)}}{4}} \Big( \frac{4\pi}{n} \Big)^{\frac{d}{2}} \det\big( (H_n^{(k)})^{-1}\big)^{1/2}
        = 2^{\frac{d}{2}}\mathrm{e}^{-\frac{Rn^{(1-2q)}}{4}}, 
    \end{align*}    
    such that $\lim_{n\to\infty} \P(A_{i,n})=0$.
    By the fact that $G_{n,N} = \bigcap_{i=1}^N A_{i,n}^c$ we get, with exploiting the independence of the events $A_{1,n},\dots,A_{N,n}$, that
    \begin{align*}
    \lim_{n \to \infty} \mathbb{P}(G_{n,N}) = \lim_{n \to \infty} \prod_{i=1}^N \mathbb{P}(A_{i,n}^c) = \prod_{i=1}^N \lim_{n \to \infty} \mathbb{P}(A_{i,n}^c)=1,
    \end{align*}
    which finishes the proof.
\end{proof}

\subsection{Proof of \cref{lem: W*}.}\label{app:lem W*}
Before we prove \cref{lem: W*} we need some auxiliary results.

\begin{lem}\label{lem:convergenceofLAweight}
    Under \cref{assum:Two} and \cref{assum:minimizer-limit} we have 
    $
        \lim_{n\to\infty} \pi_{n}^{(k)} = \pi_0(x_{*}^{(k)}).$
\end{lem}
\begin{proof}
Observe that
    \begin{align*}
        \pi_{n}^{(k)} &= \exp(-n \Phi(x_{n}^{(k)}))\, \pi_0(x_{n}^{(k)})
        = \exp(-n (\Phi(x_{n}^{(k)}) - \Phi(x_*^{(k)}) + \Phi(x_*^{(k)}))) \, \pi_0(x_{n}^{(k)}).
    \end{align*}
    Due to $\nabla^2 \Phi\big(x_*^{(k)}\big)>0$ and $\Phi \in \mathcal C^{3}(\R^d,[0,\infty))$, $\Phi$ being strongly convex in the neighborhood of $x_*^{(k)}$, by the Polyak--Lojasiewicz condition, there exists a constant $\widehat{c}>0$ such that
    \[
    \Phi(x_n^{(k)}) - \Phi(x_*^{(k)}) \leq \frac{1}{2\widehat{c}}\|\nabla \Phi(x_n^{(k)})\|^2.
    \]
    Moreover, since $\nabla \Phi(x_n^{(k)}) = \frac{1}{n} \nabla \log \pi_0(x_n^{(k)})$, we obtain
    \[
    \Phi(x_n^{(k)}) - \Phi(x_*^{(k)}) \leq \frac{1}{2\widehat{c}n^2} \|\nabla \log \pi_0(x_n^{(k)})\|^2 .
    \]
    Since $\lim_{n \to \infty} x_{n}^{(k)} = x_{*}^{(k)}$ and $\pi_0 \in \mathcal C^2(\R^d, (0,\infty))$, continuity implies $\lim_{n \to \infty} \|\nabla \log \pi_0(x_n^{(k)})\| = \|\nabla \log \pi_0(x_*^{(k)})\|$. Consequently, there exists a constant $C>0$ (independent of $n$) and $n_0\in\N$ such that for all $n \geq n_0$, we have
    \[
    \Phi(x_n^{(k)}) - \Phi(x_*^{(k)}) \leq C n^{-2}.
    \]
    Hence, since $\Phi(x_*^{(k)}) = 0$ and $\pi_0 \in \mathcal C^2(\R^d, (0,\infty))$, it follows that
    $
    \lim_{n\to\infty} \pi_{n}^{(k)} = \pi_0(x_{*}^{(k)}).
    $
\end{proof}



We provide two inequalities w.r.t. random variables that hold pointwise restricted to a suitable measurable subset of the probability space.

\begin{prop}\label{prop: upperbound_k_kappa}
For $R>0$, $q\in[0,1/2)$, $\rho,M\in(0,\infty)$ let \cref{assum:Two} and \cref{assum:minimizer-limit} be satisfied. For $N\in\mathbb{N}$ let $X_{1,n},\dots,X_{N,n}$ be iid $\lap$-distributed random variables and
$G_{n,N}$ be defined as in \eqref{eq:GnN}. Furthermore, for any $i\in\{1,\dots,N\}$ define the random index
\begin{equation}\label{eq:kappa_ni}
    \kappa_{i,n} 
    \coloneqq \argmin_{\ell\in\{1,\dots,K\}} \| X_{i,n} - x_n^{(\ell)} \|^2_{H_n^{(\ell)}}.
\end{equation}
Then, there are numbers $c_1,c_2>0$ and an $n_1\in\mathbb{N}$ such that for any $i\in\{1,\dots,N\}$ and any $n\in\mathbb{N}$ with $n\geq n_1$ we have restricted to $G_{n,N}$ that
\begin{equation}
\label{eq: ratio_estimate_restri}
\frac{\sum_{\ell=1,\,\ell\neq \kappa_{i,n}}^K \pi_n^{(\ell)}\exp\big( -\frac{n}{2} \|X_{i,n} - x_n^{(\ell)}\|^2_{H_n^{(\ell)}} \big)}{\pi_n^{(\kappa_{i,n} )}\exp\big( -\frac{n}{2} \|X_{i,n} - x_n^{(\kappa_{i,n} )}\|^2_{H_n^{(\kappa_{i,n} )}} \big)} \leq c_1 \mathrm{e}^{-c_2\,n}. 
\end{equation}
Moreover, for $\unimpw$ (defined in \eqref{weight_LM}) there are numbers $c_3 > 0$ and $n_2 \in \mathbb{N}$, such that for any $i\in\{1,\dots,N\}$ and any $n\in\mathbb{N}$ with $n\geq n_2$
we have restricted to $G_{n,N}$ that
\begin{equation}
\label{eq: weight_close_1}
 1- c_3\,n^{1-3q} \leq \unimpw (X_{i,n}) \leq 1+ c_3\,n^{1-3q}. 
\end{equation}
\end{prop}

\begin{proof}
We start with proving \eqref{eq: ratio_estimate_restri}. For this
define
\[
 \gamma :=
        \max_{j,\ell\in \{1,\dots,K\}, j\not = \ell} \frac{\pi_0(x_*^{(\ell)})}{\pi_0(x_*^{(j)})}
        \qquad \text{and} \qquad
\gamma_n 
        := \max_{j,\ell\in \{1,\dots,K\}, j\not = \ell} \frac{\pi_n^{(\ell)}}{\pi_n^{(j)}}.
\] 
Note that, by \cref{lem:convergenceofLAweight} we have $\lim_{n\to\infty} \gamma_n =   \gamma$.
Furthermore,
by \cref{assum:Two} (part \eqref{en: pairw_disj}) there is a constant $\alpha>0$, such that
    \[
        \min_{j,\ell\in\{1,\dots,K\},j\not=\ell}
        \|x_n^{(j)} - x_n^{(\ell)}\|_{H_n^{(\ell)}} 
        \geq 
        \alpha.
    \] 
   Now, choose $n_1 \in\mathbb{N}$ such that for any $n\in\mathbb{N}$ with $n\geq n_1$ holds
    $\gamma_n \leq 2\gamma$
         and 
        $
        \alpha > 4\sqrt{R} n^{-q} $.
    The latter condition is motivated by the fact that for arbitrary    
    $i\in\{1,\dots,N\}$, since 
    all random variables ($X_{i,n}$ and $\kappa_{i,n}$) are considered restricted to $G_{n,N}$, we have
    \[
        \min_{\ell \in\{1,\dots,K\}}
        \| X_{i,n} - x_n^{(\ell)} \|^2_{H_n^{(\ell)}} 
        =
        \| X_{i,n} - x_n^{(\kappa_{i,n})} \|^2_{H_n^{(\ell)}}
        \leq
        Rn^{-2q}
    \]
and, consequently, for any $\ell\in\{1,\ldots,K\}$ with $\ell\neq \kappa_{i,n}$ that
    \begin{align*}
        \|X_{i,n} - x_n^{(\ell)}\|^2_{H_n^{(\ell)}} 
        &\geq \Big| \|X_{i,n} - x_n^{(\kappa_{i,n})}\|_{H_n^{(\ell)}} - \|x_n^{(\kappa_{i,n})} - x_n^{(\ell)}\|_{H_n^{(\ell)}} \Big|^2
        \geq (\alpha - \sqrt{R}n^{-q})^2.
    \end{align*}
    Hence, for $n\geq n_1$ we obtain
    \begin{align*}
        &\frac{\sum_{\ell =1, \ell\neq \kappa_{i,n}}^K \pi_n^{(\ell)}\exp\big( -\frac{n}{2} \|X_{i,n} - x_n^{(\ell)}\|^2_{H_n^{(\ell)}} \big)}{\pi_n^{(\kappa_{i,n})}\exp\Big( -\frac{n}{2} \|X_{i,n} - x_n^{(\kappa_{i,n})}\|^2_{H_n^{(\kappa_{i,n})}} \Big)}\\
        &= \sum_{\ell =1, \ell\neq \kappa_{i,n}}^K \frac{\pi_n^{(\ell)}}{\pi_n^{(\kappa_{i,n})}}\exp\Big( -\frac{n}{2} \Big(\|X_{i,n} - x_n^{(\ell)}\|^2_{H_n^{(\ell)}} - \|X_{i,n} - x_n^{(\kappa_{i,n})}\|^2_{H_n^{(\kappa_{i,n})}} \Big) \Big)\\
        &\leq \sum_{\ell =1, \ell\neq \kappa_{i,n}}^K \gamma_n 
        \exp\Big( -\frac{n}{2} \big( ( \alpha - \sqrt{R}n^{-q})^2 - Rn^{-2q} \big)\Big)\\
        &\leq \sum_{\ell =1, \ell\neq \kappa_{i,n}}^K 2\gamma \exp\Big( -\alpha n \big( \frac{\alpha}{2} - \sqrt{R}n^{-q} \big)\Big)
        \leq  2(K-1) \gamma \exp\Big( -\frac{\alpha^2}{4} n \Big),
    \end{align*}
    such that \eqref{eq: ratio_estimate_restri} holds with $c_1=2(K-1)\gamma$ and $c_2=\alpha^2/4$.

    Now we turn to the proof of \eqref{eq: weight_close_1}. Let $i\in \{1,\dots,N\}$ be arbitrary but fixed. Now, since $\nabla I_n(x_n^{(\kappa_{i,n})}) = 0$, a second-order Taylor expansion of $I_n$ at $x_n^{(\kappa_{i,n})}$ yields
    \begin{align*}
        I_n(X_{i,n}) = I_n\big(x_n^{(\kappa_{i,n})}\big) + \frac{1}{2}\|X_{i,n}-x_n^{(\kappa_{i,n})}\|^2_{H_n^{(\kappa_{i,n})}} + R_n^{(\kappa_{i,n})}(X_{i,n}),
    \end{align*}
    where $R_n^{(\kappa_{i,n})}(X_{i,n})$ denotes the third-order remainder term. Consequently
    \begin{align*}
        \exp( -nI_n(X_{i,n})) 
        &= \exp \Big( -n \Big( I_n(x_n^{(\kappa_{i,n})}) + \frac{1}{2}\|X_{i,n}-x_n^{(\kappa_{i,n})}\|^2_{H_n^{(\kappa_{i,n})}} + R_n^{(\kappa_{i,n})}(X_{i,n}) \Big) \Big).
    \end{align*}
    As before, all random variables ($X_{i,n}$ and $\kappa_{i,n}$) are considered restricted to $G_{n,N}$, which implies $\| X_{i,n} - x_n^{(\kappa_{i,n})} \|^2_{H_n^{(\kappa_{i,n})}} \leq Rn^{-2q}$. Furthermore, since the matrices $H_n^{(k)}$ are $\rho$-uniformly positive definite and the third-order derivative tensors of $I_n$ are $M$-uniformly bounded on $\Vnq$ for all $k \in \{1,\cdots,K\}$, Taylor’s theorem implies
    \begin{equation}
    \label{eq: residual_term_Lapla}
    \big|n R_n^{(\kappa_{i,n})}(X_{i,n})\big|
    \leq 
    n \frac{6M}{3!\rho^3}\|X_{i,n}-x_{n}^{(\kappa_{i,n})}\|^3 
    \leq
    M\rho^{-3} R^{3/2} n^{1-3q}.
    \end{equation}
    Then
    \begin{align*}
        \unimpw(X_{i,n}) 
        & = \frac{\exp( -nI_n(X_{i,n}) )}{\sum_{\ell=1}^K \pi_n^{(\ell)}\exp\big( -\frac{n}{2} \|X_{i,n} - x_n^{(\ell)}\|^2_{H_n^{(\ell)}} \big)}
        = \frac{\exp \big( -n R_n^{(\kappa_{i,n})}(X_{i,n}) \big)}
        {1 + \frac{\sum_{\ell=1,\ell \neq \kappa_{i,n}}^K \pi_n^{(\ell)}\exp\big( -\frac{n}{2} \|X_{i,n} - x_n^{(\ell)}\|^2_{H_n^{(\ell)}} \big)}{\pi_n^{(\kappa_{i,n})}\exp\big( -\frac{n}{2} \|X_{i,n} - x_n^{(\kappa_{i,n})}\|^2_{H_n^{(\kappa_{i,n})}} \big)} } \\
        & \leq \exp \big( -n R^{(\kappa_{i,n})}_n(X_{i,n}) \big)
         \leq
        \exp \Big( M\rho^{-3} R^{3/2} n^{1-3q} \Big)
        \leq 1+ \exp\big(M \rho^{-3} R^{3/2}\big) n^{1-3q},
    \end{align*}
    where in the last inequality we exploited $n^{1-3q}\leq1$ (which holds by $q>1/3$) and the fact that for any $z\in [0,1]$, $a\geq 0$ holds $\mathrm{e}^{az} \leq 1+\mathrm{e}^a z$.
    Consequently
    \begin{equation}
        \label{al: lower_bnd}
      \unimpw (X_{i,n})-1 \leq  \exp\big(M \rho^{-3} R^{3/2}\big) n^{1-3q}.
    \end{equation}
     We now verify the corresponding lower bound. By \eqref{eq: ratio_estimate_restri} for $n\geq n_1$ we obtain
    \begin{align*}
        \unimpw\left(X_{i,n}\right)
        & = \frac{\exp \big( -n R_n^{(\kappa_{i,n})}(X_{i,n}) \big)}
        {1 + \frac{\sum_{\ell=1,\ell \neq \kappa_{i,n}}^K \pi_n^{(\ell)}\exp\big( -\frac{n}{2} \|X_{i,n} - x_n^{(\ell)}\|^2_{H_n^{(\ell)}} \big)}{\pi_n^{(\kappa_{i,n})}\exp\big( -\frac{n}{2} \|X_{i,n} - x_n^{(\kappa_{i,n})}\|^2_{H_n^{(\kappa_{i,n})}} \big)} }
        \geq \frac{\exp \left( -n R^{\kappa_{i,n}}_n(X_{i,n}) \right)}
        {1 + c_1 \mathrm{e}^{-c_2\,n}}.
    \end{align*}
    Exploiting \eqref{eq: residual_term_Lapla} and that for any $z\geq 0$ holds $\mathrm{e}^{-z} \geq 1-z$ we obtain
    \begin{equation*}
        \exp \left( -n R^{\kappa_{i,n}}_n(X_{i,n}) \right)
        \geq
        \exp\big(- M \rho^{-3} R^{3/2} n^{1-3q}\big)
        \geq 
        1 - M \rho^{-3} R^{3/2}\, n^{1-3q}.
    \end{equation*}
    We choose $n_2 \in \mathbb{N}$ such that $n_2\geq n_1$ and that for any $n\in\mathbb{N}$ with $n\geq n_2$ we have 
    $
    c_1 \mathrm{e}^{-c_2\,n}\leq M \rho^{-3} R^{3/2}\, n^{1-3q}.
    $
This yields   
\begin{align*}
                \unimpw\left(X_{i,n}\right)
                \geq 
                \frac{1 - M \rho^{-3} R^{3/2}\, n^{1-3q}}
        {1 + c_1 \mathrm{e}^{-c_2\,n}}
         \geq 1-M \rho^{-3} R^{3/2}\, n^{1-3q} -c_1 \mathrm{e}^{-c_2\,n}
         \geq 1-2M \rho^{-3} R^{3/2}\, n^{1-3q}.
    \end{align*}
    In combination with \eqref{al: lower_bnd}, for $c_3:=\max\{2M \rho^{-3} R^{3/2},\exp\big(M \rho^{-3} R^{3/2}\big)\}$ and $n\geq n_2$, we obtain 
    \[
    1-c_3 n^{1-3q} \leq \unimpw\left(X_{i,n}\right)
    \leq 1 + c_3 n^{1-3q}.
   \qedhere\]
\end{proof}

We collected all auxiliary tools together and are now able to prove Lemma~\ref{lem: W*}.

\begin{proof}[Proof of Lemma~\ref{lem: W*}]
Define $\overline{\unimpw}(X_{i,n}) := N^{-1} \sum_{i=1}^N \unimpw (X_{i,n})$ and note that
    \begin{align*}
        &| \ESS_n(N) - N |  = \bigg| \frac{\big(\sum_{i=1}^N \unimpw \left(X_{i,n}\right)\big)^2}{\sum_{i=1}^N \unimpw (X_{i,n})^2} - N \bigg|
        = \bigg| \frac{\big(\sum_{i=1}^N \unimpw (X_{i,n})\big)^2 - N\sum_{i=1}^N \unimpw (X_{i,n})^2}{\sum_{i=1}^N  \unimpw (X_{i,n})^2} \bigg| \\
        &= N \bigg| \frac{\big(N^{-1}\sum_{i=1}^N \unimpw (X_{i,n})\big)^2 - N^{-1}\sum_{i=1}^N \unimpw (X_{i,n})^2}{N^{-1}\sum_{i=1}^N \unimpw (X_{i,n})^2} \bigg| 
         = \frac{\sum_{i=1}^N \big( \unimpw (X_{i,n}) - \overline{ \unimpw} (X_{i,n})\big)^2}{N^{-1}\sum_{i=1}^N \unimpw (X_{i,n})^2}.
    \end{align*}
As in the proof before all random variables are considered restricted to $G_{n,N}$.
With $n_2\in\mathbb{N}$ and $c_3>0$ from Proposition~\ref{prop: upperbound_k_kappa} we have that \eqref{eq: weight_close_1} holds true. 
Consequently,  we have for $n\geq n_2$ that
\begin{align}
\label{al: again}
     1- c_3\,n^{1-3q} & \leq \unimpw (X_{i,n}) \leq 1+ c_3\,n^{1-3q},
\end{align}
which readily implies
\begin{align*}
    -2c_3 n^{1-3q} & \leq \unimpw (X_{i,n}) - \overline{\unimpw} (X_{i,n}) \leq 2c_3 n^{1-3q}.
\end{align*}
The latter inequality can be rewritten as
\[
 \vert \unimpw (X_{i,n}) - \overline{\unimpw} (X_{i,n}) \vert \leq 2c_3 n^{1-3q}.
\]
Now, we choose $n_0\in\mathbb{N}$ in such a way that $n_0 \geq n_2$ and for any $n\in\mathbb{N}$ with $n\geq n_0$ holds $1-c_3n^{1-3q}\geq 1/2$, i.e., the lower bound of \eqref{al: again} is $1/2$. Therefore, we obtain for any $n\in\mathbb{N}$ with $n\geq n_0$ that 
\begin{align*}
\sum_{i=1}^N \big( \unimpw (X_{i,n}) - \overline{ \unimpw} (X_{i,n})\big)^2& 
\leq  4Nc_3^2 n^{2-6q}
\qquad \text{and} \qquad
\frac{1}{4} \leq N^{-1}\sum_{i=1}^N \unimpw (X_{i,n})^2.
\end{align*}
Consequently, for all $n\in\mathbb{N}$ with $n\geq n_0$ with $c:=8c_3^2$ we have
$
 | \ESS_n(N) - N | \leq c Nn^{2-6q}.
$
\end{proof}

\section{Proofs regarding IMH algorithm}

\subsection{Proof of Proposition~\ref{prop:IMH_SpectralGap}}
\label{app:IMH_SpectralGap}

For any $x\in\R^d$ and Borel measurable $A\subseteq \R^d$ we have
\begin{align*}
    \mathrm{IM}_n(x,A) \geq \int_A \alpha_n(x,y) \nu({\rm d}y) = \int_A \big(\frac{1}{\unimpw(y)} \wedge \frac{1}{\unimpw(x)} \big)\, \unimpw(y) \nu({\rm d}y) 
    \geq \frac{\int_{\R^d} \unimpw(z) \nu({\rm d}z)}{\| \unimpw \|_{\infty}} \pi_n(A) 
 = \frac{\pi_n(A) }{\| \impw \|_{\infty}},
\end{align*}
where the last equality follows by $\impw(x)=\frac{\unimpw(x)}{\int_{\R^d} \unimpw(z) \nu({\rm d}z)}$. (Consequently,
the `whole state space is a small set'.) 
Setting $\gamma:=\| \impw \|_{\infty}^{-1}$ and observing that
\[
R_n(x,A) := \frac{\mathrm{IM}_n(x,A) - \gamma \pi_n(A)}{1-\gamma}
\]
is a transition kernel with stationary distribution $\pi_n$ that satisfies
\[
    \mathrm{IM}_n(x,A) = \gamma \pi_n(A) + \left(1-\gamma\right) R_n(x,A).
\]
Using that, for $f\in L^2_{0}(\post)$ with $f\neq 0$, we obtain
\begin{align*}
    &\|M_n f\|_{L^2(\pi_n)}^2
    = \int_{\mathbb R^d} \left| \int_{\mathbb R^d} f(y)\,\mathrm{IM}_n(x,{\rm d} y) \right|^2 \pi_n({\rm d}x) 
    = (1-\gamma)^2  \int_{\mathbb R^d} \left| \int_{\mathbb R^d} f(y)\,R_n(x,{\rm d}y) \right|^2 \pi_n({\rm d}x) \\
    & \leq (1-\gamma)^2 \int_{\mathbb R^d}\int_{\mathbb R^d} |f(y)|^2\,R_n(x,{\rm d}y)\,\pi_n({\rm d} x)  
    = (1-\gamma)^2 \int_{\mathbb R^d} |f(y)|^2\,\pi_n({\rm d}y) 
    = (1-\gamma)^2\|f\|_{L^2(\pi_n)}^2,
\end{align*}
where we used Jensen's inequality and stationarity of $R_n$ w.r.t. $\pi_n$. Hence,
    $\|M_n\|_{\pi_n} \leq 1-\gamma$ and the proof is finished.

\subsection{Proof of Lemma \ref{lem:EAA}}\label{app:EAA}
    By using
    $
        a \wedge b = \frac 12\left(a + b - |b-a| \right)
    $
    (which holds for all positive numbers $a,b>0$)
    we obtain 
    \begin{align*}
    \mathbb{E}_{\post\otimes\nu}[\alpha_n] & = \int_{\R^d}\int_{\R^d} \Big(1 \wedge \frac{\impw(y)}{\impw(x)} \Big) \,\nu(\mathrm dy) \post(\mathrm d x) \\
       & = \frac{1}{2}
        +
        \frac 12 \int_{\R^d}
        \int_{\R^d} \frac{\impw(y)}{\impw(x)}  \nu({\mathrm d}y)\, \impw(x) \nu( {\mathrm d}x)\\
        & \qquad -
        \frac 12 \int_{\R^d}
        \int_{\R^d} \left| 1 - \frac{\impw(y)}{\impw(x)} \right| \nu({\mathrm d}y)\, \impw(x) \nu( {\mathrm d}x)\\
        & = 1 - \frac 12 \int_{\R^d}
        \int_{\R^d} \left|\impw(x) - \impw(y) \right| \nu({\mathrm d}y) \nu( {\mathrm d}x)
        = 1 - d_{\mathrm{TV}}(\post \otimes \nu, \nu \otimes \post) 
    \end{align*}
    cf. also \cite[Theorem 1]{Neklyudov2019} for general MH algorithms.
    The first statement on the average acceptance rate then follows by
    \[
        d_{\mathrm{TV}}(\post, \nu)
        \leq
        d_{\mathrm{TV}}(\post \otimes \nu, \nu \otimes \post)
        \leq
        2d_{\mathrm{TV}}(\post, \nu)
    \]
    where the first inequality is due to the data processing inequality for general $f$-divergences and the second inequality due to applying the triangle inequality as follows
    \begin{align*}
        d_{\mathrm{TV}}(\post \otimes \nu, \nu \otimes \post)
        & \leq \frac 12 \int_{\R^d}
        \int_{\R^d}\left|\impw(x) - 1 \right| + \left|1 - \impw(y) \right| \nu({\mathrm d}y) \nu( {\mathrm d}x)
         =
        2d_{\mathrm{TV}}(\post, \nu).
    \end{align*}

    For the second statement about the lag-one autocorrelation let us first define for probability measures $\mu, \widetilde \mu$ on $\R^{2d}$ that are absolutely continuous w.r.t. the product measure $\nu_2:=\nu\otimes \nu$ with corresponding densities $\frac{{\rm d}\mu}{{\rm d}\nu_2}$, $\frac{{\rm d}\widetilde{\mu}}{{\rm d}\nu_2}$ a positive measure $\mu \wedge \widetilde \mu$ on $\R^{2d}$ by
    \[
        [\mu \wedge  \widetilde \mu](A)
        :=
        \int_A \Big(\frac{{\rm d}\mu}{{\rm d}\nu_2}(z) \wedge \frac{{\rm d}\widetilde{\mu}}{{\rm d}\nu_2}(z)\Big) \nu_2({\rm d} z),
    \]
    for arbitrary Borel measurable $A \subseteq \R^{2d}$. We have
    \begin{align*}
        \rho_n^v &= \int_{\R^d}\int_{\R^d} |v^Tx-v^Ty|^2 \Big(1 \wedge \frac{\impw(y)}{\impw(y)} \Big) \,\nu(\mathrm dy) \post(\mathrm d x)\\
        & = \int_{\R^d}\int_{\R^d} |v^Tx-v^Ty|^2 \left( \impw(x) \wedge \impw(y)\right)\nu(\mathrm dy) \nu(\mathrm d x)
        = [\pi_n\otimes \nu \wedge \nu \otimes \pi_n](|\Delta f_v|^2)
    \end{align*}
    where $f_v(x):= v^Tx$ and where for a function $f\colon \R^d \to \R$ we set
    \[
    [\pi_n\otimes \nu \wedge \nu \otimes \pi_n](|\Delta f|^2) \coloneqq \int_{\R^d}\int_{\R^d} |f(x)-f(y)|^2 \,[\pi_n\otimes \nu \wedge \nu \otimes \pi_n]({\rm d}x{\rm d}y).
    \]
    The latter notation can be readily adapted to other product measures, such that
    \[
        2\Var_{\pi_n}(f_v) = [\pi_n\otimes\pi_n \wedge \pi_n\otimes\pi_n](|\Delta f_v|^2)
    \]
    and by \eqref{eq:pif_tildepif} and Remark \ref{rem:HellBound} we have
    \begin{equation}
    \begin{aligned}
    & \left| \rho_n^v - 2\Var_{\pi_n}(f_v)\right|
     = \left|[\pi_n \otimes \nu \wedge \nu \otimes \pi_n](|\Delta f_v|^2) - [\pi_n\otimes\pi_n \wedge \pi_n\otimes\pi_n](|\Delta f_v|^2) \right|\\
    & \leq \sqrt{2[\pi_n\otimes \nu \wedge \nu\otimes \pi_n](|\Delta f_v|^4) + 2[\pi_n\otimes\pi_n \wedge \pi_n\otimes\pi_n](|\Delta f_v|^4)}\\
     & \qquad \qquad \qquad \qquad \qquad \qquad \qquad\cdot 
    d_{\mathrm H} \left(\pi_n\otimes \nu \wedge \nu\otimes \pi_n,\pi_n\otimes\pi_n \wedge \pi_n\otimes\pi_n\right).
    \end{aligned}
    \label{inequalityfromR.S.}
    \end{equation}
    Furthermore
    \begin{align*}
    & d_{\mathrm H}^2 \left(\pi_n \otimes \nu \wedge \nu \otimes \pi_n,\pi_n\otimes \pi_n \wedge \pi_n\otimes\pi_n \right)\\
    & = \int_{\R^d}\int_{\R^d} \left( \sqrt{\impw(x) \wedge \impw(y)} - \sqrt{\impw(x)\impw(y)} \right)^2 \, \nu(\mathrm dx)\nu(\mathrm dy)\\
    & \leq \int_{\R^d}\int_{\R^d} \left( \sqrt{\impw(x)} - \sqrt{\impw(x)\impw(y)} \right)^2 \, \nu(\mathrm dx)\nu(\mathrm dy) \\
    & \quad\qquad + \int_{\R^d}\int_{\R^d} \left( \sqrt{\impw(y)} - \sqrt{\impw(x)\impw(y)} \right)^2 \, \nu(\mathrm dx)\nu(\mathrm dy) 
    = 2 d_{\mathrm H}^2 \left(\pi_n,\nu \right).
    \end{align*}
    Finally, by $a\wedge b \leq a+b$ for $a,b>0$, we get
    \begin{align*}
    [\pi_n\otimes \nu \wedge \nu\otimes \pi_n](|\Delta f_v|^4)
        &= \int_{\R^d}\int_{\R^d} |v^Tx-v^Ty|^4  \left( \impw(x) \wedge \impw(y)\right)\nu(\mathrm dy) \nu(\mathrm d x) \\
        &\leq 2\int_{\R^d}\int_{\R^d}|v^Tx-v^Ty|^4 \pi_n({\rm d}x)\nu({\rm d}y)
    \end{align*}
    which by exploiting \eqref{eq:Corr_rho_v_n} implies the statement.

\bibliographystyle{plain}
\bibliography{biblio}

\end{document}